\documentclass[10pt]{amsart}
\usepackage{amsfonts}
\usepackage{amsmath}
\usepackage{amsthm}
\usepackage{amssymb}
\usepackage{latexsym}

\newtheorem*{cor}{Corollary}
\newtheorem*{lem}{Lemma}

\newtheorem*{prop}{Proposition}

\theoremstyle{definition}
\newtheorem*{defn}{Definition}
\theoremstyle{definition}
\newtheorem{thm}{Theorem}
\newtheorem*{conj}{Conjecture}
\newtheorem*{rem}{Remark}

\newenvironment{pf}{\proof}{\endproof}
\newcounter{cnt}
\newenvironment{enumerit}{\begin{list}{{\hfill\rm(\roman{cnt})\hfill}}{%
\settowidth{\labelwidth}{{\rm(iv)}}\leftmargin=\labelwidth%
\advance\leftmargin by 
\labelsep\rightmargin=0pt\usecounter{cnt}}}{\end{list}}

\theoremstyle{remark}

\numberwithin{equation}{section} 
\begin{document}

\newcommand{\thmref}[1]{Theorem~\ref{#1}}
\newcommand{\secref}[1]{Section~\ref{#1}}
\newcommand{\lemref}[1]{Lemma~\ref{#1}}
\newcommand{\propref}[1]{Proposition~\ref{#1}}
\newcommand{\corref}[1]{Corollary~\ref{#1}}
\newcommand{\remref}[1]{Remark~\ref{#1}}
\newcommand{\defref}[1]{Definition~\ref{#1}}
\newcommand{\er}[1]{(\ref{#1})}
\newcommand{\id}{\operatorname{id}}
\newcommand{\tensor}{\otimes}
\newcommand{\nc}{\newcommand}
\newcommand{\rnc}{\renewcommand}
\newcommand{\qbinom}[2]{\genfrac[]{0pt}0{#1}{#2}}
\nc{\cal}{\mathcal} \nc{\goth}{\mathfrak} \rnc{\bold}{\mathbf}
\renewcommand{\frak}{\mathfrak}
\newcommand{\desc}{\operatorname{desc}}
\newcommand{\Maj}{\operatorname{Maj}}
\renewcommand{\Bbb}{\mathbb}
\nc\bpi{{\mbox{\boldmath $\pi$}}}
\newcommand{\lie}[1]{\mathfrak{#1}}
\makeatletter
\def\section{\def\@secnumfont{\mdseries}\@startsection{section}{1}%
  \z@{.7\linespacing\@plus\linespacing}{.5\linespacing}%
  {\normalfont\scshape\centering}}
\def\subsection{\def\@secnumfont{\bfseries}\@startsection{subsection}{2}%
  {\parindent}{.5\linespacing\@plus.7\linespacing}{-.5em}%
  {\normalfont\bfseries}}
\makeatother
\def\subl#1{\subsection{}\label{#1}}
\nc{\wh}[1]{\widehat{#1}}

\nc{\Cal}{\cal} \nc{\Xp}[1]{X^+(#1)} \nc{\Xm}[1]{X^-(#1)}
\nc{\on}{\operatorname} \nc{\ch}{\mbox{ch}} \nc{\Z}{{\bold Z}}
\nc{\J}{{\cal J}} \nc{\C}{{\bold C}} \nc{\Q}{{\bold Q}}
\renewcommand{\P}{{\cal P}}
\nc{\N}{{\Bbb N}} \nc\boa{\bold a} \nc\bob{\bold b} \nc\boc{\bold
c} \nc\bod{\bold d} \nc\boe{\bold e} \nc\bof{\bold f}
\nc\bog{\bold g} \nc\boh{\bold h} \nc\boi{\bold i} \nc\boj{\bold
j} \nc\bok{\bold k} \nc\bol{\bold l} \nc\bom{\bold m}
\nc\bon{\bold n} \nc\boo{\bold o} \nc\bop{\bold p} \nc\boq{\bold
q} \nc\bor{\bold r} \nc\bos{\bold s} \nc\bou{\bold u}
\nc\bov{\bold v} \nc\bow{\bold w} \nc\boz{\bold z}

\nc\ba{\bold A} \nc\bb{\bold B} \nc\bc{\bold C} \nc\bd{\bold D}
\nc\be{\bold E} \nc\bg{\bold G} \nc\bh{\bold H} \nc\bi{\bold I}
\nc\bj{\bold J} \nc\bk{\bold K} \nc\bl{\bold L} \nc\bm{\bold M}
\nc\bn{\bold N} \nc\bo{\bold O} \nc\bp{\bold P} \nc\bq{\bold Q}
\nc\br{\bold R} \nc\bs{\bold S} \nc\bt{\bold T} \nc\bu{\bold U}
\nc\bv{\bold V} \nc\bw{\bold W} \nc\bz{\bold Z} \nc\bx{\bold X}

\title[Quantum loop modules]{\ \\[.3cm]Quantum Loop modules}\author{Vyjayanthi Chari}
\address{Department of Mathematics, University of
California, Riverside, CA 92521.}
\email{chari@math.ucr.edu}
\author{Jacob Greenstein}
\thanks{The second author was supported by a 
Marie Curie Individual Fellowship of the European Community 
under contract no.~HPMF-CT-2001-01132}
\address{Institut de Math\'ematiques de Jussieu, Universit\'e
Pierre et Marie Curie, 175 rue du Chevaleret, Plateau 7D, F-75013
Paris, France}
\email{greenste@math.jussieu.fr}
\dedicatory{Dedicated to Anthony Joseph on the occasion of his
60th birthday}
\begin{abstract}
We classify the simple infinite dimensional 
integrable modules with finite dimensional
weight spaces over the quantized enveloping algebra of an untwisted
affine algebra. We prove that these are either highest (lowest)
weight integrable modules or simple submodules of a loop module of
a finite-dimensional simple integrable module and describe the latter
class. Their characters and crystal bases theory are discussed in
a special case.
\end{abstract}
\date{June 28, 2002}
\maketitle

\section{Introduction} 
The aim of the present paper is to study
irreducible  integrable modules for quantum affine algebras,
which have finite dimensional weight spaces.
The best known examples of such representations are the highest
weight representations~$\wh{V(\lambda)}$ (cf.~\cite{Ka1,L}) 
on which the center of
the quantum affine algebra acts via a positive integer power
of~$q$. These
representations have many pleasant properties, for instance it is
known that they admit a canonical global or crystal bases. Another
family of integrable modules for the quantum affine algebra are
the finite-dimensional modules which have been studied,
amongst the others, in~\cite{AK,CP,EM,FM,FR,Ka2,N1,N2,VV}.
However, unlike the highest weight representations, these
finite-dimensional representations do not respect the natural
$\bz$-grading on the quantum affine algebra which arises from the
adjoint action of the element of the torus corresponding to the 
Euler operator. 
Thus it is natural to look for the graded analogue of the finite-dimensional 
modules. Besides, in certain cases one has to consider these infinite
dimensional modules instead of finite-dimensional ones. For example,
a finite-dimensional module cannot appear as a submodule of 
the ring of linear endomorphisms of~$\wh{V(\lambda)}$
whilst a simple integrable module with non-trivial zero weight
space can be embedded in such a ring for~$\lambda$ sufficiently
large (cf. for example~\cite{Jo1}).

Examples of infinite dimensional integrable modules which are not
highest weight modules are easy to construct. Namely, given a
finite-dimensional module $V$, one can define on the space $L(V)=
V\otimes \bc(q)[t,t^{-1}]$  in an obvious way the structure of a
graded module for the quantum affine algebra. However, even if $V$
is irreducible, the resulting  representation $L(V)$ need not remain 
so.

The irreducible finite-dimensional representations of
the quantum affine algebra are known to be parametrized by families of
polynomials in an indeterminate $u$ (cf.~\cite{CP}). For instance,
for the quantum affine
algebra associated to $\lie{sl}_2$, we have just one polynomial
$\pi(u)$. If we take $\pi(u)=(1-u)(1+u)$, and $V$ to be the
corresponding irreducible finite-dimensional representation, 
then it is not hard to see that the graded infinite dimensional module 
$L(V)$ is a direct sum of two simple components. More generally,
the module $L(V)$ is reducible (moreover,
completely reducible) if the roots of the polynomials
associated to $V$ satisfy a certain condition which involves roots
of unity. 

The paper is organized as follows. In~\secref{QL}  we show that 
if $V$ is an irreducible
finite-dimensional representation of the quantum affine algebra, then
the corresponding graded representation on
$L(V)$ is completely reducible and we describe its irreducible components
(cf.~\thmref{main1} and~\lemref{mainlemma}).
Our analysis is based on the result on irreducibility of
tensor products of finite-dimensional simple modules obtained in~\cite{C},
which allows one to construct an action of the cyclic group~$\bz/m\bz$ 
on~$L(V)$ commuting with that of the quantum affine algebra, 
the integer~$m$ being determined
by the family of polynomials corresponding to~$V$. Namely,
these turn out to be polynomials in~$u^m$.
The irreducible components of~$L(V)$, as in the classical
case (cf.~\cite{Cint, CPnew, G1}), become eigenspaces of a generator of
that group corresponding to the egeinvalues~$\zeta^i$, $i=0,\dots,m-1$, where
$\zeta$ is an~$m$th primitive root of unity. However, unlike in the
classical case when it is induced by a natural action of the
symmetric group~$S_m$ on~$V$, 
the action of~$\bz/m\bz$ is rather sophisticated and
difficult to describe explicitly.

\secref{C} is devoted to a classification of irreducible, integrable
modules with finite dimensional weight spaces of the quantum affine algebra. 
Thus our~\thmref{main2} establishes that such
a module must be isomorphic to either a highest weight module or
its (graded) dual, or to an irreducible component of the module
$L(V)$, where $V$ is some irreducible finite-dimensional module.
The corresponding result in the classical case was established 
in~\cite{Cint,CPnew}.
Besides, we show (cf.~\propref{C40}) 
that, as in the classical case (cf.~\cite{Jo2}),
all weight spaces of a simple integrable module are finite dimensional
if and only if its set of weights satisfies a boundedness condition
of~\cite{Jo1}.

The final section of the paper is concerned with the
problem of describing the characters of irreducible components of the 
module~$L(V)$ in
terms of the characters of~$V$. More precisely, we are interested
in relating the dimension of the weight spaces of an irreducible component
of $L(V)$ to the dimension of the weight spaces of $V$. In the
classical case this was done in~\cite{G1}, and the formulae
involve a modification of the Euler $\varphi$-function. 
In the quantum case we conjecture
an analogous formula and establish it in~\propref{CH65} for certain special
modules for the quantum affine algebra associated to $\lie
{sl}_{n+1}$. In order to do that, we consider a crystal $q=0$ limit of
these modules. Then the argument becomes purely
combinatorial and involves the major index of MacMahon.
In particular, we also show that the characters of these
modules in the quantum case coincide with the characters of their
classical analogues. At the end of the section we discuss a crystal
basis theory for these modules.

\section{Preliminaries}\label{P}

Throughout this paper~$\bn$ (respectively, $\bn^+$) denotes the
set of non-negative (respectively, positive) integers.
Let $q$ be an indeterminate and let $\bc(q)$ be the field of rational
functions in $q$ with complex coefficients.  For $r,m\in\bn$,
$m\ge r$, define
\begin{equation*}
[m]_q=\frac{q^m -q^{-m}}{q -q^{-1}},\qquad [m]_q!
=[m]_q[m-1]_q\ldots [2]_q[1]_q,\qquad \qbinom{m}{r}_q
= \frac{[m]_q!}{[r]_q![m-r]_q!}.
\end{equation*}

\subl{P10} Let $\frak g$ be a complex finite-dimensional simple
Lie algebra  of rank $n$ with a Cartan subalgebra~$\frak h$ and
let $W$ be the Weyl group of $\lie  g$. Set $I=\{1,2,\dots ,n\}$
and let $A=(d_ia_{ij})_{i,j\in I}$, where the~$d_i$
are positive integers, be the $n\times n$ symmetrized
Cartan matrix of $\frak g$. Let $\{\alpha_i\,:\,i\in
I\}\subset\frak h^*$ (respectively $\{\varpi_i\,:\,i\in
I\}\subset\frak h^*$) be the set of simple roots (respectively  of
fundamental weights) of $\frak g$ with respect to $\frak h$. As
usual, $Q$ (respectively, $P$) denotes the root (respectively,
weight) lattice of $\frak g$. Let~$P^+=\sum_{i\in I}\bn\varpi_i$
be the set of dominant weights and set~$Q^+=\sum_{i\in
I}\bn\alpha_i$. It is well-known that~$\lie h^*$  admits a
non-degenerate symmetric $W$-invariant bilinear form which will be
denoted by~$(\cdot\,\vert\,\cdot)$. We assume
that~$(\alpha_i\,\vert\, \alpha_i)=d_i a_{ij}$ for all~$i,j\in I$.

Let 
$$\widehat{\frak g}=\frak g\otimes\bc[t,t^{-1}]\oplus\bc c\oplus \bc d
$$ 
be the untwisted
extended affine algebra associated with~$\frak g$ and let
$\hat{A}=(d_ia_{ij})_{i,j\in\hat I}$, where~$\hat I=I\cup\{0\}$ be
the extended symmetrized Cartan matrix and $\wh W$ the affine Weyl
group. Set $\widehat{\frak h}=\frak h\oplus \bc c\oplus \bc d$.
From now on we identify $\lie h^*$ with the subspace of
$\wh{\lie h}^*$ consisting of elements which are zero on $c$ and
$d$. Define~$\delta\in \wh{\frak h}^*$ by
$$\delta(\frak h\oplus \bc c)=0,\qquad\delta(d)=1.$$
Denote by~$\theta$
the highest root of~$\lie g$ and set $\alpha_0=\delta-\theta$. Then
$\{\alpha_i\,:\,i\in\hat I\}$ is a set of simple roots for
$\wh{\lie g}$ with respect to~$\wh{\lie h}$ 
and~$\delta$ generates its imaginary roots.
The bilinear form on $\frak h^*$ extends to
a $\wh W$-invariant  bilinear form on $\wh{\lie h}^*$ which we
continue to denote by~$(\cdot\,\vert\,\cdot)$. One has $(\delta\,
\vert\,\alpha_i)=0$ and $(\alpha_i\,\vert\,\alpha_j)=d_i a_{ij}$, for
all~$i,j\in\hat I$. Define a set of fundamental
weights~$\{\omega_i\,:\, i\in \hat I\}\subset\widehat{\frak h}^*$
of~$\widehat{\frak g}$ by
the conditions~$(\omega_i|\alpha_j)=d_i \delta_{i,j}$ and
$\omega_i(d)=0$ for all~$i,j\in\hat I$. Notice that
$\varpi_i=\omega_i-\omega_0$ for all $i\in I$. Let~$\widehat P=
\sum_{i\in\hat I}\bz\omega_i\oplus\Z\delta$ (respectively,
$\widehat P^+=\sum_{i\in\hat I}\bn\omega_i\oplus\Z\delta$) be the
corresponding set of integral (respectively, dominant) weights.
Set~$P^e=P\oplus\Z\delta \subset\widehat P$. Denote by~$\widehat
Q$ the root lattice of~$\widehat{\lie g}$ and set~$\widehat
Q^+=\sum_{i\in\hat I} \bn\alpha_i$. Given $\lambda,\mu\in
\widehat{P}^+$ (respectively, $\lambda,\mu\in P^+$) we say that
$\lambda\le \mu$ if $\mu-\lambda\in \widehat{Q}^+$ (respectively,
$\mu-\lambda\in Q^+$).

\subl{P20} For $i\in\hat I$, set $q_i=q^{d_i}$ and
$[m]_i=[m]_{q^i}$. The quantum affine algebra~$\widehat\bu_q(\frak
g)$ (cf.~\cite{Be,BCP,Dr,J}) associated to $\frak g$, which will
be further denoted as~$\widehat\bu_q$, is an associative algebra
over $\bc(q)$ with generators $x_{i,r}^{{}\pm{}}\, , i\in I$,
$r\in\bz$, $K_i^{ {}\pm 1}\, ,\,i\in I$, $C^{\pm 1/2}$, $D^{\pm
1}$, $h_{i,r}\,, \,i\in I$, $r\in \bz\backslash\{0\}$, and the
following defining relations
\begin{gather*}
  \text{$C^{{}\pm{1/2}}$ are central,}\\
K_iK_i^{-1} = K_i^{-1}K_i
  =1,\quad C^{1/2}C^{-1/2} =C^{-1/2}C^{1/2} =1,\\
K_iK_j =K_jK_i,\quad
K_ih_{j,r} =h_{j,r}K_i,
\\
K_ix_{j,r}^\pm K_i^{-1} = q_i^{{}\pm a_{ij}}x_{j,r}^{{}\pm{}},
\\
  [h_{i,r},h_{j,s}]=\delta_{r,-s}\,\frac1{r}[ra_{ij}]_i\,\frac{C^r-C^{-r}}
  {q_j^{}-q_j^{-1}},
\\
[h_{i,r} , x_{j,s}^{{}\pm{}}] =
  \pm\frac1r[ra_{ij}]_i\, C^{{}\mp {|r|/2}}x_{j,r+s}^{{}\pm{}},
\end{gather*}
\begin{gather*}
  x_{i,r+1}^{{}\pm{}}x_{j,s}^{{}\pm{}} -q_i^{{}\pm
    a_{ij}}x_{j,s}^{{}\pm{}}x_{i,r+1}^{{}\pm{}}=q_i^{{}\pm
    a_{ij}}x_{i,r}^{{}\pm{}}x_{j,s+1}^{{}\pm{}}
  -x_{j,s+1}^{{}\pm{}}x_{i,r}^{{}\pm{}},
\\
[x_{i,r}^+ ,
  x_{j,s}^-]=\delta_{i,j}\,\frac{ C^{(r-s)/2}\,\psi_{i,r+s}^+ -
    C^{-(r-s)/2}\, \psi_{i,r+s}^-}{q_i^{} - q_i^{-1}},
\\
\sum_{\pi\in\Sigma_m}\sum_{k=0}^m(-1)^k\qbinom{m}{k}_{i}
  x_{i, r_{\pi(1)}}^{{}\pm{}}\ldots x_{i,r_{\pi(k)}}^{{}\pm{}}
  x_{j,s}^{{}\pm{}} x_{i, r_{\pi(k+1)}}^{{}\pm{}}\ldots
  x_{i,r_{\pi(m)}}^{{}\pm{}} =0,\qquad\text{if $i\ne j$},
\end{gather*}
for all sequences of integers $r_1,\ldots, r_m$, where $m
=1-a_{ij}$, $\Sigma_m$ is the symmetric group on $m$ letters, and
the $\psi_{i,r}^{{}\pm{}}$ are determined by equating powers of
$u$ in the formal power series
$$
\sum_{r=0}^{\infty}\psi_{i,\pm r}^{{}\pm{}}u^{{}\pm r} = K_i^{{}\pm 1}
\exp\left(\pm(q_i^{}-q_i^{-1})\sum_{s=1}^{\infty}h_{i,\pm s} u^{{}\pm s}
\right).
$$
The algebra $\wh\bu_q({\lie g})$ is $\bz$-graded, with the
$l^{th}$-graded piece being,
$$(\wh\bu_q({\lie g}))_l=\{x\in\wh\bu_q(\lie g)\,:\, DxD^{-1}=q^l x\}.$$
The subalgebra of $\widehat\bu_q(\lie g)$ generated by the
elements $x_{i,0}^\pm$, $i\in I$ is isomorphic to the quantized
enveloping algebra $\bu_q(\lie g)$ of $\lie g$. Let
$\wh\bu_q^\pm(\gg)$ (respectively $\wh\bu_q^\pm(\ll)$, $\wh\bu^\pm_q(0)$)
be the subalgebra of $\wh\bu_q$ generated by the elements
$x^+_{i,r}$, $i\in I$, $\pm r\in\bn$ (respectively, $x^-_{i,r}$, $i\in I$,
$\pm r\in\bn$, $h_{i,r}$, $i\in I$, $\pm r\in\bn^+$). Let
$\wh\bu_q^\circ$ be the subalgebra generated by $K_i^{\pm 1}$,
$i\in I$, $D^{\pm1}$ and~$C^{\pm1/2}$. We will need the following
result which was established in~\cite{BCP}.
\begin{prop}\label{triangle} The subspaces
$\wh\bu_q^\pm=\wh\bu_q^\pm(\ll)\wh\bu_q^\pm(0)\wh\bu_q^\pm(\gg)$ are
subalgebras of $\widehat\bu_q$ and 
$$\wh\bu_q =
\wh\bu_q^-\wh\bu_q^\circ\wh\bu_q^+.
$$
\end{prop}

\subl{P25}
We will also need another presentation of~$\widehat\bu_q$. Namely,
after~\cite{Be,J}, the algebra $\widehat\bu_q$ is isomorphic to an
associative $\bc(q)$-algebra generated by $E_i, F_i, K_i^{{}\pm
1}\,:\,i\in\hat I$, $D^{\pm 1}$ and central elements~$C^{\pm1/2}$
satisfying the following relations:
\begin{gather*}
\text{$C=K_0\displaystyle\prod_{i\in I} K_i^{r_i}$, where~$\theta=
\displaystyle\sum_{i\in I} r_i\alpha_i$, $r_i\in\N^+$}
\\
K_iK_i^{-1}=K_i^{-1}K_i=1,\qquad K_iK_j=K_jK_i,\\
  DD^{-1}=D^{-1}D=1,\qquad K_iD=DK_i,\\
K_iE_j K_i^{-1}=q_i^{ a_{ij}}E_j,\qquad
K_iF_j K_i^{-1}=q_i^{-a_{ij}}F_j,\\
DE_j D^{-1}=q^{ \delta_{j0}}E_j,\qquad
DF_j D^{-1}=q^{-\delta_{j0}}F_j,\\
[E_i, F_j]=\delta_{ij}\,\frac{K_i-K_i^{-1}}{q_i^{}-q_i^{-1}},
\end{gather*}
\begin{gather*}
  \sum_{r=0}^{1-a_{ij}}(-1)^r\qbinom{1-a_{ij}}{r}_i
(E_i)^rE_{j}(E_{i})^{1-a_{ij}-r}=0\
  \qquad \text{if $i\ne j$},\\
\sum_{r=0}^{1-a_{ij}}(-1)^r\qbinom{1-a_{ij}}{r}_i
(F_{i})^rF_{j}(F_{i})^{1-a_{ij}-r}=0\
  \qquad \text{if $i\ne j$}.
\end{gather*}
The element ~$E_i$ (respectively, $F_i$), $i\in I$ corresponds
to~$x_{i,0}^+$ (respectively, $x_{i,0}^-$). In particular, the
$E_i,F_i,K_i^{\pm1}\,:\,i\in I$ generate~$\bu_q(\frak g)$.

It is well-known that $\widehat\bu_q$ is a Hopf algebra over
$\bc(q)$ with  the co-multiplication being given in terms of
generators $E_i,F_i,K_i^{\pm1} \,:\,i\in\hat I$ by the following
formulae
$$
\Delta(E_i)=E_i\tensor K_i^{-1}+1\tensor E_i,\qquad
\Delta(F_i)=F_i\tensor 1+K_i\tensor F_i,
$$
the~$K_i^{\pm1}$ being group-like. Although explicit formulae for
the co-multiplication on generators ~$x_{i,r}^\pm$, $h_{i,r}$ are
not known, we will not need these in the present paper and so we
say no more about it.
\begin{lem}\label{comul}
Take~$x\in (\widehat\bu_q)_k$ and write~$\Delta(x)=x_1\tensor x_2$
in the summation notation. Then we may assume
that~$x_i\in(\widehat\bu_q)_{k_i}$ where~$k_1+k_2=k$.
\end{lem}
\begin{pf}
The assertion is obvious for the generators~$E_i,F_i,K_i^{\pm1}$.
Since $\Delta$ is an algebra homomorphism, it holds for any
polynomial in these generators, which is homogeneous with respect
to~$D$.
\end{pf}

\subl{P30} Let $\bu_q^e$ be the extended quantum loop algebra,
namely the graded quotient of $\widehat\bu_q$ by the graded ideal
generated by~$C^{\pm 1/2}-1$. The Hopf algebra structure on
$\wh\bu_q$ descends to a Hopf algebra structure on $\bu_q^e$. Let
$\bu_q$ be the $\bc(q)$-subalgebra of~$\bu_q^e$ generated by the
elements $x^{\pm}_{i,k}$, $h_{i,r}$, $K_i, K_i^{-1}$, $i\in I$,
$k,r\in\bz$, $r\ne 0$. It is easy to see that $\bu_q$ is a Hopf
subalgebra of $\bu_q^e$. Let $\bu_q^e(0)$ be the subalgebra of
$\bu_q^e$ generated by elements $h_{i,r}$, $K_i^{\pm 1}$, $i\in
I$, $r\in\bz$ and $D^{\pm 1}$. Clearly $\bu^e_q(0)$ is a
$\bz$-graded subalgebra of $\bu^e_q$. Let $\bu_q(>)$ (respectively, 
$\bu_q(<)$, $\bu_q(0)$) be the subalgebra of~$\bu_q^e(0)$
generated by the elements $x^+_{i,r}$, $i\in I$, $r\in\bz$ (respectively,
$x^-_{i,r}$, $i\in I$, $r\in\bz$, $h_{i,r}$, $i\in I$, $r\in\bz\setminus
\{0\}$. Then~(cf. for example~\cite{CP})
$$\bu_q^e=\bu_q(<)\bu^e_q(0)\bu_q(>).$$

For $i\in I$,  set
\begin{equation*}
h^\pm_i(u)=\sum_{k=1}^\infty \frac{q^{\pm k}h_{i,\pm
k}}{[k]_i}u^k,
\end{equation*}
and
$$P^\pm_i(u)=\exp(- h^\pm_i(u))=
\exp\Big(-\sum_{k=1}^\infty \frac{q^{\pm k}h_{i,\pm
k}}{[k]_i}u^k\Big).
$$
Let $P_{i,\pm r}$ be the coefficient of
$u^r$ in $P_i^\pm (u)$. One can show (cf.~\cite{BCP}) that
the~$P_{i,r}\,:\,i\in I$, $r\in\Z$ generate $\bu_q(0)$. Moreover, by
\cite{BCP} monomials in the $P_{i,r}$, $i\in I$, $r\in\bz$ (or
equivalently, monomials in the $h_{i,r}$, $i\in I$, $r\in\bz$)
form a basis of $\bu_q(0)$. 
\begin{lem}{\label{chi}}
Let~$\chi:\bu_q(0)\to\bc(q)[t,t^{-1}]$ be a non-trivial
homomorphism of $\bz$-graded algebras.
\begin{enumerit}
\item
There exists a unique~$m>0$ such
that the image of~$\chi$ equals~$\bc(q)[t^m,t^{-m}]$.
\item
Suppose that the image of~$\chi$
equals~$\bc(q)[t^m,t^{-m}]$. Then there exist $i_0\in I$ such that
$\chi(P_{i_0,\pm m})\not=0$ and the kernel of $\chi$
is generated by the~$P_{i,r}$, $i\in I$, $r\not=0\pmod m$ and
by the elements of the form
$$
P_{i,\pm sm}-(\chi(P_{i_0,\pm m})^{-s}\chi(P_{i,\pm sm}))
P^s_{i_0,\pm m},\qquad i\in I, s\in\bn^+.
$$
\end{enumerit}
\end{lem}
\begin{proof}
Let $m,n\in\bn^+$ be minimal such that both~$t^m$, $t^{-n}$ lie in
the image of~$\chi$. Then~$t^{m-n}\in\operatorname{Im}\chi$. Since
either $0\le m-n < m$ or $0\le n-m < n$ it follows that~$m=n$ which proves 
the
first part. For the second part, it is enough to observe that the
elements listed in the assertion are homogeneous and lie in the kernel
of~$\operatorname{Im}\chi$.
\end{proof}

\subl{P50}
In the rest of this section we summarize some general results from
the representation theory of $\widehat\bu_q$ and $\bu_q(\frak g)$, which
will be used later.

For any $\widehat\bu_q$-module $\widehat{V}$ and any
$\mu=\sum_{i\in\hat I}\mu_i\omega_i+l\delta\in\widehat{P}$, set
\begin{equation*}
\widehat{V}_{\mu}=\{ v\in \widehat{V}\,:\, D.v=q^{l}v,\, K_i.v
=q_i^{\mu_i}v,\,\forall\, i\in \hat{I}\}.
\end{equation*}
If~$\widehat V_\mu\not=0$ we say that~$\mu$ is a weight of~$\widehat V$.
The set of weights of~$\widehat V$ will be denoted by~$\Omega(\widehat V)$.
The module $\widehat V$ is said to be an
{\em admissible  module of type~$1$\/} if
\begin{equation*}
\widehat{V}=\bigoplus_{\mu\in \widehat{P}}
\widehat{V}_{\mu},
\end{equation*}
and $\dim \widehat{V}_\mu<\infty$ for all $\mu\in\widehat{P}$.
One has analogous definitions of admissible modules of type~$1$ for the
algebras $\bu_q^e$ (with~$\widehat P$ replaced by~$P^e$),
$\bu_q$ and $\bu_q(\frak g)$ (with~$\widehat P$ replaced by~$P$).
From now on, all modules will be assumed to be of type~$1$.
A $\widehat\bu_q$-module $\widehat{V}$ is
integrable if for all $i\in \hat{I}$ the elements $E_i$ and $F_i$
act locally nilpotently on $\widehat{V}$. Similarly, one can define
integrable $\bu_q^e$, $\bu_q$ and~$\bu_q(\lie g)$-modules.

\subl{P60}
We now recall the construction
of highest weight integrable modules over $\widehat\bu_q$
and~$\bu_q(\frak g)$. We  work with the presentation of~$\widehat\bu_q$
described in~\ref{P30}.

The following result can be found in \cite[3.5]{L}.
\begin{prop}
For every $\lambda=\sum_{i\in\hat I}\lambda_i \omega_i+k\delta\in
\widehat P^+$ {\em({\em respectively $\lambda=\sum_{i\in I}
\lambda_i\varpi_i\in P^+$})} there exists a unique, up to an
isomorphism, simple integrable $\widehat\bu_q$-module
$\widehat V(\lambda)$ {\em({\em respectively $\bu_q(\frak g)$-module
$V(\lambda)$})} of type~$1$ which is generated by an element
$v_\lambda$ satisfying
$$
E_i.v_\lambda=0,\qquad K_i.v_\lambda=q_i^{\lambda_i}v_\lambda,\qquad
D.v_\lambda=q^k v_\lambda,\qquad
F_i^{\lambda_i+1}.v_\lambda=0,\qquad \forall\,i\in\hat I
$$
{\em({respectively, $E_i.v_\lambda=0$,
$K_i.v_\lambda=q_i^{\lambda_i}v_\lambda$,
$F_i^{\lambda_i+1}.v_\lambda=0$, $\forall\,i\in I$})}.
\end{prop}

\subl{P70} The modules  $\widehat V(\lambda)$ and $V(\lambda)$ are
the quantum analogues of the corresponding modules $\widehat
V_{c\ell}(\lambda)$ and $V_{c\ell}(\lambda)$ over, respectively,
$\widehat{\frak g}$ and $\frak g$, whose characters are given by
the Weyl-Kac formula. Namely (cf.~\cite[Theorem~4.12]{L1}), for
all $\nu\in \widehat P$, (respectively $\nu\in P$), we have $\dim
\widehat V(\lambda)_\nu =\dim \widehat V_{c\ell}(\lambda)_\nu$
(respectively, $\dim V(\lambda)_\nu=\dim V_{c\ell}(\lambda)_\nu$).
In particular, both $\widehat V(\lambda)$ and $V(\lambda)$ have
finite-dimensional weight spaces. Furthermore, a standard argument
from the representation theory of~$\lie g$ yields the following
\begin{lem}\label{weight} Let $\mu\in P^+$ be such that $\lambda-\mu\in 
Q^+$.
Then $\dim V(\lambda)_\mu\ne 0$.
\end{lem}
Finally, we will need the following
\begin{prop}[{cf.~\cite[3.5~and~6.3]{L}}]\label{lus}
\begin{enumerit}
\item Assume that $V=\bigoplus_{\nu\in\widehat P}
V_\nu$ is an integrable admissible  $\widehat\bu_q$-module
such that~$\Omega(V)$ is contained in the set
$\bigcup_{i=1}^k\{\mu_i-\nu\,:\,\nu\in\widehat Q^+\}$, for
some~$k\in\N^+$ and for some~$\mu_1,\dots,\mu_k\in\widehat P^+$.
Then
$$
V \cong \bigoplus_{\lambda\in\widehat P^+} m(\lambda) \widehat
V(\lambda),
$$
for some non-negative integers $m(\lambda)$. Furthermore, as a
$\bu_q(\lie g)$-module, $\widehat V(\lambda)$ is a direct sum of
simple finite dimensional highest weight modules~$V(\mu)$ for
various~$\mu\in P^+$.

\item Let~$M$ be an
integrable $\bu_q(\lie g)$-module. Then~$M$ is a sum of
simple~$\bu_q(\lie g)$-modules of form~$V(\lambda)$ for various
$\lambda\in P^+$.
\end{enumerit}
\end{prop}
\begin{cor} \label{intfin}
Let~$M$ be an integrable~$\bu_q(\lie g)$-module with finite dimensional
weight spaces. Then~$M$ is finite-dimensional.
\end{cor}
\begin{pf}
Let~$r\in\bn^+$ be minimal such that~$r\varpi_i\in Q^+$ for all~$i\in I$.
Then any~$\lambda\in P^+$ can be written as~$\lambda'+\eta$ where~$\eta\in
Q^+$ and~$\lambda'\in P^+_r:=\{\nu=\sum_{i\in I}k_i\varpi_i\in P^+\,:\,
0\le k_i<r\}$. Evidently, $P^+_r$ is a finite set. On the other hand,
by~\lemref{weight}, for all~$\lambda\in P^+$, there exists~$\lambda'\in
\Omega(V(\lambda))\cap P^+_r$. Since~$\dim M_{\lambda'}<\infty$ and
the sum of~$V(\lambda)$ is direct, it follows that the multiplicity
of each~$V(\lambda)$ in~$M$ is finite and that the set of~$\lambda\in P^+$
such that~$V(\lambda)$ occurs in~$M$ is also finite. It remains to
apply~(ii) of the above Proposition.
\end{pf}

\subl{P80} Since ~$\widehat\bu_q$ is a Hopf algebra,  given a
$\widehat\bu_q$-module~$V=\bigoplus_{\nu\in\widehat P} V_\nu$, we
can endow~$V^*$ with a structure of a~$\widehat\bu_q$ module via
the antipode. If~$\dim V_\nu<\infty$ for all~$\nu\in\widehat P$
then~$V^\#=\bigoplus_{\nu\in\widehat P} V^*_\nu$, is a
$\widehat\bu_q$-submodule of~$V^*$. The module~$V^\#$ is called
the graded dual of~$V$. One can prove that the graded dual
of~$\widehat V(\lambda) \,,\,\lambda\in\widehat P^+$ is the unique
simple integrable module generated by an element~$v^*_\lambda$
such that
\begin{equation*}\label{liw} F_i.v^*_\lambda=0,\qquad
K_i.v^*_\lambda=q_i^{-\lambda_i}v^*_\lambda, \qquad
D.v^*_\lambda=q^{-k} v^*_\lambda,\qquad E_i^{\lambda_i+1}
v^*_\lambda=0.
\end{equation*}
Clearly, $\widehat V(\lambda)^\#$ is an integrable module. Results
analogous to the ones above hold for the modules $\widehat
V(\lambda)^\#$. Finally, note that the element $C$ acts on
$\widehat V(\lambda)$ as $q^r\id$  and on $\widehat V(\lambda)^\#$
as $q^{-r}\id$ where~$r=(\lambda\,\vert\,\delta)$. Notice
that~$r>0$ unless~$\lambda\in\wh P^+\cap P^e=\bz\delta$. In the
latter case both~$\wh V(\lambda)$ and~$\wh V(\lambda)^\#$ are
one-dimensional.

\subl{P90}
Let~$M$ be an integrable~$\bu_q(\lie g)$ or $\wh\bu_q$-module.
Following~\cite[5.2]{L}, one
can define $\bc(q)$-linear endomorphisms $T_i\,:\,i\in I$ (respectively,
$i\in\hat I$) of~$M$ satisfying~$T_i M_\lambda=M_{s_i\lambda}$.
Moreover, by~\cite[Chapter~39]{L}, the~$T_i$ satisfy the
relations of the braid group associated with~$W$ (respectively, $\wh W$).
In particular, the set of weights of~$M$ is $W$ (respectively~$\wh W$)
invariant. Moreover, if~$M$ is admissible, then we have $\dim M_\lambda=
\dim M_{w\lambda}$ for all~$w\in W$ (respectively, $w\in\wh W$).

\section{Quantum loop modules}\label{QL}

In this section we study
a family of irreducible integrable modules of $\widehat\bu_q$ on
which $C$ acts as the identity and hence these are actually modules for the
extended  quantum loop algebra $\bu_q^e$.

\subl{QL10}
Given a~$\bu_q$ module $V$ of type~1, one can easily verify that
the following formulae  define a structure of a $\bu_q^e$-module
of type~1, which we denote by~$L(V;d)$, on the vector space
$L(V)=V\otimes\bc(q)[t,t^{-1}]$. Namely, fix~$d\in\Z$ and define
for all~$k,r\in\bz$, $x\in(\bu_q^e)_k$ and~$v\in V$
$$
x(v\otimes t^r)=(xv)\otimes t^{r+k},\qquad D(v\otimes t^r)=q^{d+r}\,
v\otimes t^r.
$$
We set  $L(V)=L(V;0)$ and call it  the quantum loop module
associated to $V$. 
\begin{lem}\label{trivlem}
Let $V$ be a cyclic $\bu_q$ module generated by an
element $v\in V$. Then $L(V)$ is generated as a $\bu_q^e$-module
by the elements $v\otimes t^r$, $r\in\bz$. 
\end{lem}
\begin{pf}
Immediate.
\end{pf}
If $V$ is a finite-dimensional $\bu_q$-module then the
corresponding  loop module $L(V)$ is an integrable
$\bu^e_q$-module.

The main result of this section is the
following
\begin{thm}\label{main1} Let $V$ be an irreducible
finite-dimensional $\bu_q$-module and let $d\in\bc$. Then, there
exists $v\in V$ and a unique $m\in\bn^+$ such that as
$\bu_q^e$-modules we have, $$
L(V;d)=\bigoplus_{s=0}^{m-1}\bu_q^e.(v\otimes t^s), $$ where
$\bu^e_q.(v\otimes t^s)$ is an irreducible $\bu_q^e$-module for
all $0\le s\le m-1$.
\end{thm}
We prove this theorem in the remainder of this section.
For simplicity of notation,
we assume that
$d=0$, the general case being identical.

\subl{QL20}
We need several results about irreducible
finite-dimensional representations of $\bu_q$ which we now
recall (cf.~\cite{CP}).
Let $$\cal{A}=\{\pi=\sum_{m\ge
0}\pi_mu^m\in\bc(q)[[u]]\,:\,\pi(0)=1\}.$$
\begin{defn}\label{hw} We say that a $\bu_q$-module $V$ is $\ell$-highest
weight, with highest weight $(\lambda, \bpi^\pm)$, where
$\lambda=\sum_{i\in I}\lambda_i\varpi_i$,
$\bpi^\pm=(\pi_1^\pm(u),\cdots,\
\pi_n^\pm(u))\in\cal{A}^n$,  if there
exists $0\ne v\in V_\lambda$ such that $V=\bu_q.v$ and
\begin{equation*}
x_{i,k}^+. v= 0,\qquad K_i.v =q_i^{\lambda_i}v,\qquad
P^\pm_i(u).v=\pi_i^\pm(u)v,
\end{equation*}
for all $i\in I$,
$k\in\bz$. Such an element $v$ is called a highest weight vector.
\end{defn}
If $V$ is an $\ell$-highest weight module,
then in fact $V=\bu_q(<).v$ and so
\begin{equation*} V_\mu\ne 0\implies \mu=\lambda-\eta\quad(\eta\in Q^+).
\end{equation*}
For any  $\lambda\in P^+$, $\bpi^\pm\in\cal{A}^n$, there exists
a unique (up to an isomorphism) irreducible highest weight
$\bu_q$-module  with highest weight $(\lambda,\bpi^\pm)$. Write
$\pi_i^\pm(u)=\sum_{r\ge0} \pi_{i,r}^\pm u^r$ and
let
$I(\lambda,\bpi^\pm)$ be the left ideal in $\bu_q$ generated by
$\bu_q(>)_+$ and the elements $P_{i,\pm r}-\pi^\pm_{i,r}$, 
$K_i-q_i^{\lambda_i}\,
:\,i\in I$, $r\in\bn^+$ and let
$M(\lambda,\bpi^\pm)$ be the quotient of $\bu_q$
by $I(\lambda,\bpi^\pm)$. Let~$\bar v_\bpi$ be the
canonical image of~$1\in\bu_q$
in~$M(\lambda,\bpi^\pm)$.
This module is $\ell$-highest weight
and has a unique irreducible quotient which we denote as
$V(\lambda,\bpi^\pm)$. Let~$v_\bpi$ be the canonical image of~$\bar v_\bpi$
in~$V(\lambda,\bpi^\pm)$.
\begin{lem}\label{phi}
Let $\bpi^\pm\in\cal{A}$. Take $a\in\bc(q)^\times$ and set
$\bpi^\pm_a =(\pi^m_1(au),\cdots ,\pi^\pm_n(au))$. Then there
exists a canonical isomorphism $\phi_{\bpi^\pm,
a}:V(\lambda,\bpi^\pm)\to V(\lambda,\bpi^\pm_a)$ of
$\bu_q(\mathfrak g)$-modules, such that
\begin{enumerit}
\item For all $x\in(\bu_q)_k$, $v\in V(\lambda,\bpi^\pm)$ we
have,
$$
\phi_{\bpi^\pm, a}(xv)=a^{k}x\phi_{\bpi^\pm, a}(v).
$$
\item For all $a,b\in\bc(q)^\times$,
$$\phi_{\bpi_a^\pm, b}\circ \phi_{\bpi^\pm,a}=\phi_{\bpi^\pm, ab}.$$
\end{enumerit}
\end{lem}
\begin{pf} Given $a\in\bc(q)^\times$, let
$\phi_a$ be the
graded algebra automorphism of~$\bu_q^e$ defined on generators by
$$
\phi_a(x^\pm_{i,k})=
a^{k}x^\pm_{i,k},\qquad \phi_a(h_{i,k})=a^{k} h_{i,k},\qquad
\phi_a(K_i)=K_i,\qquad \phi_a(D)=D,
$$
for all~$i\in I$ and~$k\in\bz$.
Evidently, $\phi_a^{-1}=\phi_{a^{-1}}$. It follows immediately
from the definitions that $\phi_a$ maps $I(\lambda,\bpi^\pm)$
to $I(\lambda,\bpi_a^\pm)$ and hence induces an isomorphism  of
vector spaces $\phi_a: M(\lambda,\bpi^\pm)\to
M(\lambda,\bpi^\pm_a)$. Since $\phi_a(x)=x$ if $x\in(\bu_q)_0$ it
follows that this map is an isomorphism of $(\bu_q)_0$-modules,
in particular, a map of $\bu_q(\frak g)$-modules. Let
$L(\lambda,\bpi^\pm)$ be the maximal $\bu_q$-submodule of
$M(\lambda,\bpi^\pm)$.
We claim that~$\phi_a(L(\lambda,\bpi^\pm))$ is a~$\bu_q$-submodule
of~$L(\lambda,\bpi^\pm_a)$, even though~$\phi_a$ is not a~$\bu_q$-module
map. Indeed, take~$y\in\bu_q$ such that~$y+I(\lambda,\bpi^\pm)\in
L(\lambda,\bpi^\pm)$. Then for all~$x\in\bu_q$,
\begin{align*}
x\phi_a(y+I(\lambda,\bpi^\pm))&=x\phi_a(y)+I(\lambda,\bpi^\pm_a)
=\phi_a(\phi_a^{-1}(x)y)+ I(\lambda,\bpi^\pm_a)\\
&=\phi_a(\phi_a^{-1}(x)y+ I(\lambda,\bpi^\pm))\in
\phi_a(L(\lambda,\bpi^\pm)).
\end{align*}
Since $\phi_a(\bar v_{\bpi})=\bar v_{\bpi_a}$, we conclude that
$\phi_a$ factors through to a map $V(\lambda,\bpi^\pm)\to
V(\lambda,\bpi^\pm_a)$, which is the desired isomorphism~$\phi_{\bpi^\pm,a}$.
The properties (i) and (ii) of $\phi_{\bpi^\pm,a}$ are
immediate from the fact the algebra automorphism $\phi_a$
satisfies these conditions.
\end{pf}

\subl{QL30}
The following was proved in \cite{CP} (see also \cite{CPweyl}).
\begin{thm}\label{classify}
Assume that  $\lambda\in P^+$ and $\bpi^\pm\in
\cal{A}^n$ satisfy the following
\begin{enumerit}
\item $\bpi^+=(\pi_1,\cdots ,\pi_n)$ is an
$n$-tuple of polynomials,
\item $\lambda_i= \deg
\pi_i$,
\item $\pi_i^-(u) =
u^{\deg\pi_i}\pi_i(u^{-1})/
\left.\left(u^{\deg\pi_i}\pi_i(u^{-1}\right))\right|_{u=0}$.
\end{enumerit}
Then $V(\lambda,\bpi^\pm)$ is an irreducible finite-dimensional
$\ell$-highest weight $\bu_q$-module with highest weight 
$(\lambda,\bpi^\pm)$.
Moreover, these exhaust the irreducible finite-dimensional
$\ell$-highest weight modules.
\end{thm}

\noindent
\hskip\parindent From now, we shall only consider  $n$-tuples of polynomials
$\bpi=(\pi_1,\cdots, \pi_n)$, which have constant term~$1$ and
split over $\bc(q)$ and we denote by~$V(\bpi)$ the
irreducible  $\ell$-highest weight module corresponding to~$\bpi$.

We also need the following (cf.~\cite{CPweyl})
\begin{prop}
\label{tensor} Assume that $V(\bpi)$ and $V(\bpi')$ are
irreducible
highest weight $\bu_q$-modules and also that $V(\bpi)\otimes
V(\bpi')$ is irreducible. Then $$V(\bpi)\otimes
V(\bpi')\cong
V(\bpi\bpi')\cong V(\bpi')\otimes V(\bpi),$$
where~$\bpi\bpi'=(\pi_1^{}\pi'_1,\dots,\pi_n^{}\pi'_n)$.
\end{prop}

\subl{QL35}
Given such an $n$-tuple $\bpi$, and $d\in\bz$,
define a homomorphism
$\chi_{\bpi,d}:\bu_q^e(0)\to\bc(q)[t,t^{-1}]$ by extending the
assignment,
\begin{equation}\label{chipi}\chi_{\bpi,d}(P_{i,\pm r})=
\pi^\pm_{i,r}t^{\pm r}, \qquad \chi_{\bpi,d}(K_i) =q_i^{\deg\pi_i}, \qquad
\chi_{\bpi,d}(D)=q^d,
\end{equation}
for all $i\in I$, $r\in\bz$ to a homomorphism of $\bz$-graded
algebras. We set $\chi_\bpi=\chi_{\bpi,0}$.
\begin{lem}\label{factor}
Let~$\bpi=(\pi_1,\dots,\pi_n)$ be an $n$-tuple of polynomials
satisfying~$\pi_i(0)=1$.
Suppose that $\chi_\bpi:\bu_q^e\to\bc(q)[t^m,t^{-m}]\to 0$,
for some $m>1$. Then for all
$i\in I$, $\pi_i$ is actually a  polynomial in $u^m$. In
particular, there exists an $n$-tuple of polynomials
$\bpi^0=(\pi_1^0,\cdots ,\pi_n^0)$ such that
$\pi_i^0(0)=1$ and
$$\pi_i=\prod_{s=0}^{m-1}\pi_i^0(\zeta^su),$$
where~$\zeta\in\bc$ is an $m^{\rm th}$ primitive root of unity.
\end{lem}
\begin{proof}
By the definition of~$\chi_{\bpi}$ and~\lemref{chi}, $\pi^\pm_{i,r}=0$
unless $m$ divides~$r$. It follows that~the $\pi$ are polynomials in~$u^m$.
In particular, $m$ divides~$\deg \pi_i$, hence we can write~$\deg\pi_i=
m_i m$ for some~$m_i\in\bn$.
Let~$a_i$ be a root of~$\pi_i$ in the algebraic closure of~$\bc(q)$.
Then, evidently, $\zeta^s a_i$ is a root of~$\pi_i$ for all~$s=0,\dots,m-1$.
Since~$\pi_i(0)=1$, $a_i\not=0$, and so we can write
$$
\pi_i(u)=\prod_{r=1}^{m_i} \prod_{s=0}^{m-1} (1-b_{i,r}\zeta^s u),
$$
for some~$b_{i,r}\in\overline{\bc(q)}$. Set
$$
\pi_i^0(u)=\prod_{r=1}^{m_i} (1-b_{i,r}u).
$$
It remains to observe that~$\pi_i^0\in\bc(q)[u]$.
\end{proof}
Given $\bpi$ such that $\chi_\bpi:\bu_q^e\to\bc(q)[t^m,t^{-m}]\to
0$ for some $m>1$, let $\bpi^0$ be the $n$-tuple of polynomials
defined in the previous Lemma.

\subl{QL40}
Our further analysis is based on the following result
which was established in~\cite[Theorem~3]{C}.
\begin{thm}{\label{tenspr}} Let $\{\bpi_j=(\pi_{j,1},\cdots ,\pi_{j,n})\,:
\,1\le j\le k\}$ be a set of $n$-tuples of
polynomials with constant term one  which are split over $\bc(q)$.
Assume further, that for all $1\le j,j'\le k$ and $l\in I$,
we have $$a\ne bq^{\bz},$$ where $a$ and $b$ are arbitrary roots
of $\pi_{j,l}$ and $\pi_{j',l}$. Then, the tensor product
$V(\bpi_1)\otimes\cdots \otimes V(\bpi_k)$ is an irreducible
$\bu_q$-module.
\end{thm}
\begin{prop}\label{isomorphism}
Let $V(\bpi)$ be an irreducible finite-dimensional representation
of $\bu_q$, and assume that
$\chi_\bpi:\bu_q^e\to\bc(q)[t^m,t^{-m}]\to 0$, for some $m>1$.
Then there exists an $n$-tuple of polynomials~$\bpi^0$ such that
$$
V(\bpi)\cong V(\bpi^0(u))\otimes V(\bpi^0(\zeta
u))\otimes\cdots\otimes V(\bpi^0(\zeta^{m-1}u)).
$$
\end{prop}
\begin{proof}
This follows immediately from~\thmref{tenspr} and~\lemref{factor}.
\end{proof}

\subl{QL50}
It follows from~\thmref{tenspr} and~\propref{isomorphism} that there exists 
a
non-trivial isomorphism of $\bu_q$-modules,
\begin{multline*}
\tau_{\bpi^0}: V(\bpi^0(u))\otimes V(\bpi^0(\zeta
u))\otimes\cdots\otimes V(\bpi^0(\zeta^{m-1}u))\to\\
V(\bpi^0(\zeta u))
\otimes V(\bpi^0(\zeta^{2}u))\otimes\cdots\otimes V(\bpi^0(u)),
\end{multline*}
which maps the
tensor product of highest weight vectors on the left to   the
corresponding element on the right. Set
$$
\eta_{\bpi^0}=
(\phi_{\bpi^0_\zeta,\zeta^{m-1}}\otimes
\phi_{\bpi^0_{\zeta^2},\zeta^{m-1}}\otimes\cdots\otimes
\phi_{\bpi^0,{\zeta^{m-1}}}
)\circ\tau_{\bpi^0}.
$$
Then $\eta_{\bpi^0}$ is a $\bu_q(\frak
g)$-module endomorphism of $V(\bpi)\cong V(\bpi^0(u))\otimes V(\bpi^0(\zeta
u))\otimes\cdots\otimes V(\bpi^0(\zeta^{m-1}u))$ and maps a highest weight
vector to its multiple. We may assume, without loss of generality,
that~$\eta_{\bpi^0}$ maps a highest weight vector to itself.
\begin{lem}\label{eta1}
Suppose that~$x\in(\bu_q)_k$ for some $k\in\Z$. Then
$$
\eta_{\bpi^0}(x.v)=\zeta^{-k} x.\eta_{\bpi^0}(v).
$$
In particular, $\eta_{\bpi^0}^m=\id$.
\end{lem}
\begin{pf}
Let~$\Delta$ denote the standard co-multiplication on~$\bu_q$ 
(cf.~\ref{P30})
and set
$\Delta^r=(\Delta\tensor\id^{\tensor r-1})\circ\cdots
\circ (\Delta\tensor\id)\Delta:\bu_q\to\bu_q^{\tensor r+1}$.
Suppose that~$x\in (\bu_q)_k$. Then we can write,
using the summation notation of Sweedler and~\lemref{comul}, 
$\Delta^{m-1}(x)=x_1\tensor\cdots
\tensor x_m$ where~$x_j\in (\bu_q)_{k_j}$, $k_j\in\bz$ 
and~$k_1+\cdots+k_m=k$.
Therefore,
\begin{align*}
&\eta_{\bpi^0}(x v)=(\phi_{\bpi^0_\zeta,\zeta^{m-1}}\otimes
\phi_{\bpi^0_{\zeta^2},\zeta^{m-1}}\otimes\cdots\otimes
\phi_{\bpi^0,{\zeta^{m-1}}}
)(x_1\tensor\cdots\tensor x_m)\tau_{\bpi^0}(v)\\
&=\zeta^{(m-1)k} (x_1\tensor\cdots\tensor x_m)
(\phi_{\bpi^0_\zeta,\zeta^{m-1}}\otimes
\phi_{\bpi^0_{\zeta^2},\zeta^{m-1}}\otimes\cdots\otimes
\phi_{\bpi^0,{\zeta^{m-1}}})\circ \tau_{\bpi^0}(v)\\
&=\zeta^{-k} x.\eta_{\bpi^0}(v).
\end{align*}
Since~$V(\bpi)$ is simple, it follows that~$V(\bpi)=\bu_q v_\bpi$.
On the other hand, $\eta_{\bpi^0}v_\bpi=v_\bpi$, hence by the above
$V(\bpi)$ is a direct sum of eigenspaces of~$\eta_{\bpi^0}$
and all the eigenvalues of~$\eta_{\bpi^0}$ are $m$th roots of
unity. If follows immediately that~$\eta_{\bpi^0}^m=\id$.
\end{pf}

\subl{QL55}
Extend  of~$\eta_{\bpi^0}$ to~$L(V(\bpi))$ by
setting~$\eta_{\bpi^0}(v\tensor t^r)=\zeta^r\eta_{\bpi^0}(v) \tensor t^r$.
From now on we will denote this map by~$\eta$.
\begin{lem}\label{eta2}
The map~$\eta:L(V(\bpi))\to L(V(\bpi))$ is a $\bu_q^e$~module endomorphism.
\end{lem}
\begin{pf}
Take~$x\in(\bu_q^e)_k$ for some~$k\in\bz$.
Then, for all~$v\in V(\bpi)$ and~$r\in\Z$
\begin{alignat*}{2}
\eta(x(v\tensor t^r))&=\eta( xv\tensor t^{r+k})=
\zeta^{r+k}\eta_{\bpi^0}(xv)\tensor t^{r+k}\qquad&\\
&=\zeta^{r} x.\eta_{\bpi^0}(v)\tensor t^{r+k}&\text{by~\lemref{eta1}}\\
&=x.\eta(v\tensor t^r).
\end{alignat*}
Since every~$x\in\bu_q^e$ can be written, uniquely, as a finite sum
of homogeneous elements, we conclude that~$\eta$ commutes with the
action of~$\bu_q^e$.
\end{pf}

\subl{QL60}
Our~\thmref{main1} is an immediate consequence of the following
\begin{lem}\label{mainlemma}
Let~$L^s(V(\bpi))\subset L(V(\bpi))$ be the eigenspace of~$\eta$
corresponding to the eigenvalue~$\zeta^s$, $s=0,\dots,m-1$. Then
$$
L(V(\bpi))=\bigoplus_{s=0}^{m-1}L^s(V(\bpi))
=\bigoplus_{s=0}^{m-1}\bu_q^e.(v_\bpi\otimes
t^s),
$$
where the summands are simple~$\bu_q^e$-modules.
\end{lem}
\begin{pf} 
Since~$\eta_{\bpi^0}(v_\bpi)=v_\bpi$, it follows that
$\eta(v_\bpi\tensor t^s)=\zeta^s v_\bpi\tensor t^s$, hence $v_\bpi\otimes
t^s\in L^s(V(\bpi))$.
Furthermore, $\eta\in\operatorname{End}_{\bu_q^e}L(V(\bpi))$
whence~$L^s(V(\bpi))$ is a $\bu_q^e$-submodule  of~$L(V(\bpi))$.
Finally, observe that, by~\lemref{trivlem},
every proper submodule of~$L(V(\bpi))$
meets $v_{\bpi}\tensor\bc(q)[t,t^{-1}]$. Therefore,
$\bu_q^e(v_{\bpi}\tensor t^s)$ is simple.
\end{pf}

\begin{rem}
Although we have worked over $\bc(q)$,
the results of the section go through if we specialize $q$ to be a
complex number which is {\em not} a root of unity.
\end{rem}

\section{Classification of irreducible  integrable  modules for
$\widehat\bu_q$.}\label{C}

\subl{C10}
We begin with the following definition which is analogous to~\defref{hw}.
\begin{defn} A $\bu^e_q$-module~$V$ of type~$1$ is called
$\ell$-highest weight  if there exists a vector $v\in V$,
$\lambda=\sum_{i\in I}\lambda_i\varpi_i+ d\delta\in
P^e$ and a maximal graded ideal
$\mathfrak M$ in $\bu_q^e(0)$ such that $$x^+_{i,k}.v =0,\qquad
\mathfrak M.v =0,\qquad K_i.v =q_i^{\lambda_i}v,\qquad Dv=q^{d}v $$
for all~$i\in I$ and~$k\in\bz$.
\end{defn}
It is not hard to see that the set of  maximal graded ideals in
$\bu_q(0)$ is in bijective correspondence with the set of graded
ring homomorphisms $\chi:\bu_q(0)\to\bc(q)[t,t^{-1}]$. Given such
a $\chi$ and $\lambda\in P^e$, one can define in the
obvious way a universal highest weight module $M(\lambda, \chi)$.
Namely $M(\lambda,\chi)$  is  the left $\bu_q^e$-module obtained
by taking the quotient of $\bu_q^e$ by the left ideal generated by
the elements $x^+_{i,r}$, $i\in I,r\in\bz$, $\ker\chi$,
$K_i-q_i^{\lambda_i}$, $i\in I$ and $D-q^{d}$. Let~$\bar
v_{\lambda,\chi}$ be the image of~$1\in\bu_q^e$
in~$M(\lambda,\chi)$. Standard methods show that $M(\lambda,\chi)$
has a unique irreducible quotient $V(\lambda,\chi)$. Denote the
canonical image of~$\bar v_{\lambda,\chi}$ in $V(\lambda,\chi)$
by~$v_{\lambda,\chi}$.

We need the following result which
was proved in~\cite[Proposition~3.5]{CPqa}.
\begin{prop}\label{gi}  For
all $s\ge 0$, we  have, $$(x^+_{i,\pm 1})^s(x^-_{i,0})^s =P_{i,\pm
s}+\,\text{\rm terms in $\bu_q^e\bu_q(>)_+$}.$$\end{prop}

\begin{thm}\label{lclass} The $\bu_q^e$-module $V(\lambda,\chi)$ is
integrable if and only if there exists an
$n$-tuple of polynomials $\bpi=(\pi_1,\cdots ,\pi_n)$ in an
indeterminate $u$ with constant term one, and $d\in\bc$,  such
that
$$
\lambda_i=\deg\pi_i,\qquad \chi=\chi_{\bpi,d},
$$
where~$\chi_{\bpi,d}$ was defined in~\eqref{chipi}.
\end{thm}
\begin{pf} Suppose that $V(\lambda,\chi)$ is integrable. Then we
have $$ (x^-_{i,0})^{s}.v_{\lambda,\chi}=0,\qquad i\in I,\quad
s\ge \lambda_i+1, $$ 
whence 
$$(x^+_{i,\pm
1})^{s}(x^-_{i,0})^{s}.v_{\lambda,\chi}=0,\qquad i\in I,\quad s\ge
\lambda_i+1. $$ 
Using~\propref{gi}, we conclude that 
$$
P_{i,s}.v_{\lambda,\chi}=0,\qquad i\in I,\qquad |s|\ge
\lambda_i+1, 
$$ 
that is, $P_{i,s}\in\ker\chi$ for all $i\in I$ and
$|s|\ge \lambda_i+1$. Furthermore, since
$(x^-_{i,0})^{\lambda_i}.v_{\lambda,\chi} \ne 0$, and
$x^-_{i,-1}(x^-_{i,0})^{\lambda_i}.v_{\lambda,\chi} =0$, it
follows from the representation theory of $\bu_q(\lie{sl}_2)$ that
$$(x^+_{i,1})^{\lambda_i}(x^-_{i,0})^{\lambda_i}.v_{\lambda,\chi}\ne 0,$$
hence $P_{i,\lambda_i}\ne 0$. Thus, we can define an $n$-tuple
of polynomials $\bpi$ by $$\pi_i(u)=\left(\frac
{1}{\chi(P_{i,\lambda_i})t^{-\lambda_i}}\right)\sum_{r\ge 0}
(\chi(P_{i,r})t^{-r})u^r. $$ One can prove then as in \cite{CPqa},
that the eigenvalues of $P_{i,-r}$ are given by the $n$-tuple
$\bpi^-$ defined in the previous section. It is now easy to check
that $\chi=\chi_{\bpi,d}$ where~$\chi(D)=q^d$ and the result follows.

For the converse statement it suffices to prove that given an
$n$-tuple of polynomials with constant term one and $d\in \bz$,
there exists an irreducible integrable module with highest
weight $\chi_{\bpi,d}$. Consider the module $L(V(\bpi);d)$ defined
in~\ref{QL10}. It is clear that the element $v_\bpi\otimes 1$
generates a highest weight $\bu_q^e$-module which is integrable.
Further, by~\thmref{main1}, it follows that this module is
irreducible and hence is isomorphic to $V(\deg\bpi, \chi_{\bpi,d})$.
\end{pf}
We note the following consequence.
\begin{cor}\label{class} Any irreducible
integrable $\ell$-highest  weight $\bu_q^e$-module is isomorphic
to a simple submodule of a quantum loop module, and in particular
has finite-dimensional weight spaces.
\end{cor}

\subl{C20} {\em We are aiming to prove a classification result
which  similar to the one obtained in~\cite{Cint,CPnew}. This
can be done only if we work over an algebraically closed field
containing~$\bc(q)$. The classification result also
holds if  we specialize $q$ to be a non-zero complex number which
is not a root of unity. We will assume without further comment
that we are in one of these situations until the end of this
section.}

We begin with the following
\begin{prop} \label{ghw}
Let $V=\bigoplus_{\nu\in \widehat{P}} V_\nu$ be an integrable
$\widehat\bu_q$-module of type~$1$, such that
$\dim V_\nu<\infty$ for all $\nu\in\widehat P$. Then there exists
$\lambda\in \widehat{P}$ such that $V_\lambda\ne 0$ and
$V_{\lambda+\eta}=0$ for all $\eta\in Q^+\setminus\{0\}$. In particular,
$(\lambda\,\vert\,\alpha_i)\ge0$ for all~$i\in I$.
\end{prop}
\begin{pf} Suppose that for each~$\mu\in\Omega(V)$ there exists~$\nu\in
Q^+$ such that~$\mu+\nu\in\Omega(V)$. Fix some~$\mu\in\Omega(V)$. Then
there exists an infinite sequence $\{\eta_r\}_{r\ge 1}$ such
that~$\eta_r \le \eta_{r+1}$ and $\mu+\eta_r\in\Omega(V)$ for all~$r\ge1$.
Set $W_r:=\bu_q(\frak g)V_{\mu+\eta_r}$. Then~$W_r$ is an integrable
$\bu_q(\frak g)$-module with finite-dimensional weight spaces
and hence by~{\propref{lus}(ii)} and~\corref{lus}
$$
W_r\cong \bigoplus_{\mu_{r,s}\in P^+}m(\mu_{r,s})V(\mu_{r,s}),
$$
where~$m(\mu_{r,s})\in\bn$ and are non-zero for finitely
many~$\mu_{r,s}$. Choose $s_1$ such that $\nu_1:=\mu_{1,s_1}>\mu$
and $m(\mu_{1,s_1})\ne 0$. Such~$s_1$ exists since~$\mu+\eta_1$ is
a weight of~$W_1$. Furthermore, let $r_2$ be the smallest positive
integer so that there exists $s_2$ with
$\nu_2:=\mu_{r_2,s_2}>\mu$, $\mu_{1,s_1}\ne \mu_{r_2,s_2}$ and
$m(\mu_{r_2,s_2})\ne 0$. Notice that $r_2$ always exists since the
module $W_1$ is finite-dimensional and the maximal weights which
occur in $W_r$ keep increasing. Repeating this process, we obtain
an infinite collection of elements $\nu_k>\mu$, $k\ge 1$ such that
such that $V(\nu_k)$ is isomorphic to an irreducible $\bu_q(\frak
g)$-submodule $W(\nu_k)$ of $V$. By~\lemref{weight} it follows
that $V_{\mu}\cap W(\nu_k)\ne 0$ for all~$k\ge1$. Since all
the~$\nu_k$ are distinct, the sum of~$W(\nu_k)$ is direct, which
contradicts the finite-dimensionality of~$V_\mu$. In
particular~$\lambda+\alpha_i$ is not a weight for all~$i\in I$. It
follows that~$V_\lambda$ generates a highest weight $\bu_q(\frak
g)$-module, whence~$\lambda$ is dominant with respect to the
$\alpha_i$, $i\in I$.
\end{pf}

\subl{C25} Let $\bu^+_{i,s}$ (respectively, $\bu^-_{i,s}$), $i\in
I$, $s\in\bz$ be the subalgebra of~$\wh\bu_q$ generated by
$E=x^+_{i,s}$, $F=x^-_{i,-s}$ and~$K=C^s K_i$ (respectively, by
$E=x^-_{i,s}$, $F=x^+_{i,-s}$ and~$K=C^s K_i^{-1}$). One can
easily see from the defining relations (cf.~\ref{P20}) that
$\bu^\pm_{i,s}$ is isomorphic to~$\bu_q(\frak{sl}_2)$, $E$, $F$
and~$K$ being the standard generators of the latter algebra.

The following simple Lemma will be used repeatedly in the proof of
our main theorem.
\begin{lem}\label{sl2}
Let~$V$ be an integrable admissible~$\wh\bu_q$-module.
Assume that $\mu\in\Omega(V)$
satisfies~$(\mu\,\vert\,s\delta\mp\alpha_i)>0$ for some~$i\in I$ 
and~$s\in\bz$.
Then $\mu\pm\alpha_i-s\delta\in\Omega(V)$.
\end{lem}
\begin{pf}
Recall that if $M$ is a finite dimensional~$\bu_q(\frak{sl}_2)$
module and~$Km=q^k m$ for some~$m\in M$, with $k>0$, then
$Fm\not=0$. Consider a $\bu_q(\lie{sl}_2)$-module
$\bu^\mp_{i,s}V_\mu$. Since $C^s K_i^{\mp1}$ acts on~$V_\mu$
by~$q^k$ with~$k=(\mu\,|\,s\delta\mp\alpha_i)$ positive, it
follows that~$x^\pm_{i,-s}V_\mu\not=0$. Since
~$x^\pm_{i,-s}V_\mu\subset V_{\mu\pm\alpha_i-s\delta}$ the result
follows.
\end{pf}

\subl{C30}
Now we will prove the main result of this section.
\begin{thm}\label{main2}
Let $V=\bigoplus_{\nu\in\widehat P}V_\nu$ be an irreducible
integrable $\widehat\bu_q$-module of type~$1$ and assume that
$\dim V_\nu<\infty$ for all $\nu\in\widehat P$
and let~$r\in\bz$ be
such that~$Cv=q^rv$ for all $v\in V$. Suppose that~$\dim V>1$.
\begin{enumerit}
\item If $r>0$, then there exists $\lambda\in\widehat P^+$ such
that $V\cong V(\lambda)$.
\item If $r=0$ then, there exists an
$n$-tuple of polynomials $\bpi$ with constant term~$1$
and~$d\in\bz$ such that $V$ is isomorphic to an irreducible
component of $L(V(\bpi);d)$.
\item If $r<0$,  then there exists
$\lambda\in\widehat P^+$ such that $V\cong V(\lambda)^\#$.
\end{enumerit}
\end{thm}

\begin{pf}
By~\propref{ghw} we can choose
$\lambda\in \widehat
P\in\Omega(V)$ such that $\lambda+\eta$ is not a weight of~$V$ for
all $\eta\in Q^+$ and~$(\lambda\,\vert\,\alpha_i)\ge0$ for all~$i\in I$.
Given~$\eta=\sum_{i\in I} k_i\alpha_i\in Q^+$, set~$\on{ht}\eta:=\sum_{i\in 
I}
k_i$.

(i) Assume that~$r>0$. Then there exists~$m\in\bn$ such
that~$V_{\lambda+s\delta}=0$ for all~$s>m$. Indeed, otherwise we
can choose $s>0$ such that $\lambda+s\delta\in\Omega(V)$ and
$(\lambda+s\delta\,|\,s\delta-\alpha_i)=
rs-(\lambda\,|\,\alpha_i)>0$ for all $i\in I$. It follows
from~\lemref{sl2} that $\lambda+\alpha_i$ is a weight of~$V$ for
all~$i\in I$, which contradicts the choice of~$\lambda$. Next, we
prove that the $\bu_q(\lie g)$-module
$M=\wh\bu_q^+(\ll)\wh\bu_q^+(0)\wh\bu_q^+(\gg)V_{\lambda+m\delta}$
is finite-dimensional. Since
$$\wh\bu_q^+(0)\wh\bu_q^+(\gg)\subset\wh\bu_q^+(\gg)\wh\bu_q^+(0),$$
it suffices, by the choice of $m$, to prove that
$\wh\bu_q^+(\ll)\wh\bu_q^+(\gg)V_{\lambda+m\delta}$ is
finite-dimen\-sional.

Let us prove first that
\begin{equation}\label{>>}\wh\bu_q^+(\gg)V_{\lambda+m\delta}=
\bu_q^+(\lie g)V_{\lambda+m\delta},
\end{equation}
where~$\bu_q^+(\lie g)$ is the subalgebra of~$\bu_q(\lie g)$
generated by the~$x_{i,0}^+\,:\,i\in I$. First, suppose
that~$x^+_{i,s}V_{\lambda+m\delta}\not=0$ for some~$i\in I$ and~$s\in\bn^+$.
Then~$\lambda+(m+s)\delta+\alpha_i$ is a weight of~$V$.
Since~$(\lambda+(m+s)\delta+\alpha_i\,|\,\alpha_i)=
(\lambda\,|\,\alpha_i)+2 d_i>0$, it follows from~\lemref{sl2}
that~$\lambda+(m+s)\delta\in\Omega(V)$ which is a contradiction.
Then using induction on~$\on{ht}\eta$ we conclude that
$\lambda+\eta+(m+s)\delta$, $\eta\in Q^+$ is not a weight of~$V$
for all~$s>0$, whence
$$\wh\bu_q^+(\gg)V_{\lambda+m\delta}=\bu_q^+(\frak
g)V_{\lambda+m\delta}\subset \bu_q(\lie g)V_{\lambda+m\delta}.$$
Yet $\bu_q(\lie g)V_{\lambda+m\delta}$ is
finite-dimensional by~\corref{intfin} and~\eqref{>>} is proved. 

In particular, there
exists a finite set $\{\eta_k\}\subset Q^+$ such
that~$\bu_q^+(\frak g)V_{\lambda+m\delta}\subset \bigoplus_k
V_{\mu_k}$ where~$\mu_k= \lambda+\eta_k+m\delta$ and that
$M\subset\bigoplus_{k}\wh\bu_q^+(\ll)V_{\mu_k}$. Since $\dim
\wh\bu_q^+(\ll)_{-\eta+s\delta}<\infty$ 
for all~$\eta\in Q^+$ and~$s\in\bn$, to prove that $M$ is
finite-dimensional, it is now sufficient to prove that
\begin{equation}\label{<<}
\wh\bu_q^+(\ll)_{-\eta+s\delta}V_{\mu_k}=0 \end{equation} for all but
finitely many $\eta\in Q^+$ and $s\in\bn^+$, which is an immediate 
consequence of the following two assertions:

\noindent
$1^\circ.$ The set~$\{\eta\in Q^+\,:\, \mu_k-\eta+l\delta\in
\Omega(M)\,\text{for some $l>0$}\}$ is finite.

\noindent
$2^\circ.$ For every $\eta\in Q^+$,
$M_{\mu_k-\eta+l\delta}= 0$ for~$l$ sufficiently large.

In order to prove $1^\circ$, assume it to be false and
observe that the weights of $M$ are all of the form
$\mu_k-\eta+l\delta\,:\,\eta\in Q^+$. Since $M$ is an integrable
$\bu_q(\frak g)$-module, it follows (cf.~\ref{P90}) that
$\Omega(M)$ is $W$-invariant. Since 
the set of the~$\mu_k$ is finite, as is the group~$W$, and~$P$, $\delta$ are
preserved by~$W$, we can
always choose $\eta\in Q^+$ such that~$\mu_k-\eta+l\delta$ is a
weight and the $W$-orbit of~$\mu_k-\eta+ l\delta$ contains an
element which is not of the form $\mu_r-\eta'+l\delta$ for
some~$r$ and~$\eta'\in Q^+$, which is a contradiction. 

To prove $2^\circ$ we proceed
by induction on $\on{ht}\eta$. If $\on{ht}\eta=1$, that
is~$\eta=\alpha_j$ for some~$j\in I$, and
$\mu_k-\eta+l\delta\in\Omega(M)$ for infinitely many $l$, choose
$l$ large enough so that~$\mu_k-\alpha_j+l\delta\in \Omega(M)$
and~$(\mu_k-\eta+l\delta\,|\,(l-1)\delta-\alpha_j)=
r(l-1)-(\mu_k\,|\,\alpha_j)+2d_j>0$. Applying~\lemref{sl2} we
conclude that $\mu_k+\delta$ is a weight of~$M$, which is a
contradiction.
For the inductive step, suppose that $\eta\in Q^+$
with~$\on{ht}\eta>1$. Then there exists~$j\in I$ such that~$\eta-\alpha_j\in
Q^+$. Furthermore, by the induction hypothesis there exists~$N$ such
that~$\mu_k-(\eta-\alpha_i)+l\delta\notin\Omega(M)$ for all~$l>N$. 
Choose~$s>0$
so that~$\mu_k-\eta+l\delta\in\Omega(M)$ whilst~$
(\mu_k-\eta+l\delta\,|\,(s+1)\delta-\alpha_j)=r(s+1)-(\mu_k\,|\,\alpha_j)+
(\eta\,|\,\alpha_j)>0$, where~$l=N+s+1$.
Then~\lemref{sl2} yields
$\mu_k-(\eta-\alpha_j)+(N+1)\delta\in\Omega(M)$, which
is a contradiction by the choice of~$N$.

Since~$V$ is simple, it follows from Proposition~\ref{triangle}
that $\Omega(V)\subset \Omega(M)-\wh Q^+$. Further, any $\mu\in
\wh P$ with $(\mu\,|\,\delta)>0$ is $\wh W$ conjugate to  an element in
$\wh P^+$. Since~$V$ is integrable, $\Omega(V)$ is~$\wh W$-invariant 
hence $\Omega(V)$ is contained in a finite union of
cones of the form $\mu-\wh Q^+$ where $\mu\in\wh P^+$. 
In order to complete the proof of~{(i)}, it only remains 
to apply~\propref{lus}.

(ii)  If ~$r=0$ then $\Omega(V)$ is contained in~$P^e$. We prove
first that there exists $\mu\in \Omega(V)$ such that
$x_{i,m}^+V_\mu=0$ for all $i\in I$, $m\in\bz$.

Suppose that $x^+_{i,m}V_\lambda\ne 0$ for some~$i\in I$
and~$m\in\bz$. Set $\mu=\lambda+\alpha_i+m\delta$. Suppose further
that $x^+_{j,s}.V_\mu\ne 0$ for some $j\in I$, $s\in\bz$. Observe
that, if~$i\not=j\in I$ then either~$a_{ij}+2>0$ or~$a_{ji}+2>0$
and obviously $a_{ii}+2>0$. Thus we may assume, without loss of
generality, that~$a_{ij}+2>0$. Then
$(\lambda+\alpha_i+\alpha_j\,|\,\alpha_i)=(\lambda\,|\,\alpha_i)
+d_i(a_{ij}+2)>0$ and so by~\lemref{sl2} we conclude that
$\lambda+\alpha_j\in\Omega(V)$, which is a contradiction by the
choice of $\lambda$.

Thus, we have proved that $V=\widehat\bu_q v_0$ where $v_0$
satisfies
$$
x^+_{i,r}.v_0=0,\qquad\forall\,i\in I,\quad r\in\bz.
$$
Since $\bu_q^e=\bu_q(<)\bu^e_q(0)\bu_q(>)$ it follows from
standard arguments that $\bu_q^e(0) v_0$ must be an irreducible
$\bu_q^e(0)$-module.  Since $\bu_q^e(0)$ is a $\bz$-graded
commutative algebra, it follows that the irreducible graded
representations must just be the quotient of $\bu_q^e(0)$ by a
maximal graded ideal~$\frak M$ of~$\bu_q^e(0)$, which annihilates~$v_0$.
This proves that $V$ is a
$\ell$-highest weight module of $\bu_q^e$ and~{(ii)} now follows
from Corollary \ref{class}.

(iii) This case is similar to the first one.
\end{pf}

\subl{C40} Assume that the bilinear form on~$\wh{\lie h}^*$ is
normalized in such a way that its values on~$\widehat P$ are
rational (for, it is sufficient to assign a rational value
to~$(\omega_0\,|\,\omega_0)$). Let~$M=\bigoplus_{\nu \in\wh P}
M_\nu$ be a $\wh\bu_q$-module of type~$1$. We say,
following~\cite{Jo1,Jo2}, that~$M$ is a bounded module
if~$(\nu\,\vert\,\nu)\le N$ for some~$N\in\bz$ and for
all~$\nu\in\Omega(M)$. In particular, if~$M$ is simple, then this
upper bound is attained (cf.~\cite[7.2]{Jo1}), that is there
exists~$\lambda\in \Omega(M)$ such
that~$(\nu\,\vert\,\nu)\le(\lambda\,\vert\, \lambda)$ for
all~$\nu\in\Omega(M)$. We call such a~$\lambda$ maximal.

Observe that a bounded module is necessarily integrable. Indeed,
if both~$\mu$ and~$\mu+n\alpha_i$ are weights of~$M$,
then~$(\mu+n\alpha_i\,\vert\,\mu+n\alpha_i)=
(\mu\,\vert\,\mu)+2n(\mu\,\vert\,\alpha_i)+n^2(\alpha_i\,\vert\,\alpha_i)
\le N$ which imposes a bound on~$|n|$.
It is shown in~\cite{Jo2} that a simple bounded $\wh{\lie g}$-module
is admissible, being necessarily of one of the types described
in~\cite{Cint,CPnew}.

We will now establish a similar result for~$\wh\bu_q$-modules.
\begin{prop}
Let~$V$ be a simple bounded $\wh\bu_q$-module and suppose that~$\lambda$ is
its maximal weight and that~$\dim V\not=1$.
\begin{enumerit}
\item If~$(\lambda\,\vert\,\delta)>0$, then~$V\cong V(\mu)$, where~$\mu
\in\wh W\lambda\cap\wh P^+$.
\item If~$(\lambda\,\vert\,\delta)=0$, then~$V$ is isomorphic to a
simple submodule of~$L(V(\bpi);d)$ for some $n$-tuple~$\bpi$
of polynomials with constant term~$1$ and~$d\in\bz$.
\item If~$(\lambda\,\vert\,\delta)<0$, then~$V$ is isomorphic 
to~$V(\mu)^\#$,
where~$\mu\in \wh W(-\lambda)\cap \wh P^+$.
\end{enumerit}
In particular, a simple bounded $\wh\bu_q$-module is admissible.
\end{prop}
\begin{pf}
Since the form~$(\cdot\,\vert\,\cdot)$ is~$\wh W$-invariant and~$V$
is integrable, any element of the~$\wh W$-orbit of~$\lambda$ is also a maximal
weight of~$V$.

(i) Suppose that~$(\lambda\,\vert\,\delta)>0$. Then
the~$W$-orbit of~$\lambda$ contains~$\mu\in\wh P^+$. Since~$\mu$ is
maximal by the above remark, it follows that~$\mu+\alpha_i$ is not
a weight of~$V$ for all~$i\in\hat I$. Indeed, $(\mu+\alpha_i\,\vert\,
\mu+\alpha_i)=
(\mu\,\vert\,\mu)+2(\mu\,\vert\,\alpha_i)+(\alpha_i\,\vert\,\alpha_i)>
(\mu\,\vert\,\mu)$, 
which  contradicts the maximality of~$\mu$. Therefore, $E_iV_\mu=0$
for all~$i\in\hat I$, that is~$V$ is a highest weight module
of highest weight~$\mu$. Since~$V$ is simple, the claim follows 
immediately from~\propref{lus}.

(ii) Suppose that~$(\lambda\,\vert\,\delta)=0$. Then $\lambda\in P^e$ and 
so~$\wh W\lambda\cap\wh P^+$ is empty. However, since~$W$ is a finite group,
we can always conjugate~$\lambda$ to some $\mu\in P^e$ such that~$(\mu\,|\,
\alpha_i)\ge0$ for all~$i\in I$. 
Since~$(\mu+\alpha_i+s\delta\,\vert\,\mu+\alpha_i+s\delta)=(\mu\,\vert\,\mu)+
2(\mu\,\vert\,\alpha_i)
+(\alpha_i\,\vert\,\alpha_i)>(\mu\,\vert\,\mu)$ and~$\mu$ is a maximal 
weight, we 
conclude that~$\mu+\alpha_i+s\delta$ 
is not a weight of~$V$ for all~$i\in I$ and~$s\in\bz$. It follows
that~$x^+_{i,s}V_\mu=0$ for all~$i\in I$ and~$s\in\bz$. The rest
of the argument repeats that of the proof of the second part of~\thmref{main2}.

(iii) The argument repeats that of~(i).
\end{pf}

\section{Characters of quantum loop modules in $\lie{sl}_{n+1}$ case}\label{CH}

\subl{CH10}
In the classical case, we have the following result (cf.~\cite{G1}). Let
$V=\bigotimes_{i=1}^k V(\lambda_i)^{\tensor m_i m}$ be a finite
dimensional~$L\frak g=\frak g\tensor\C[t,t^{-1}]$-module
with evaluation parameters~$a_{i,r}\zeta^s$
where~$\lambda_i\in P^+$ and are distinct, $a_{i,r}\in\C$, $r=1,\dots,m_i$,
$s=0,\dots,m-1$ and~$a_{i,r}/a_{i',r'}$ is not an~$m$th root of unity.
Let~$L(V)$ be the corresponding loop $\widehat{\frak g}$-module. Then
the formula of~\cite[Theorem~4.4]{G1} yields
\begin{equation}\label{CH10.1}
\dim L^s(V)_{\nu+r\delta}=\frac1m\sum_{d|m}
\varphi_{r-s}(d) \dim V^{1/d}_{\nu/d},
\end{equation}
where~$V^{1/d}:=\bigotimes_{i=1}^k V(\lambda_i)^{\tensor m_im/d}$
and
$$
\varphi_k(d)=\varphi(d)\,\frac{\mu(d/\gcd(d,k))}{\varphi(d/\gcd(d,k))},
$$
where~$\varphi$ is the Euler function and~$\mu$ is the classical
M\"obius function, $\mu(k)=0$ if~$k$ is divisible by a square
and~$\mu(k)=(-1)^l$ if~$k$ is a product of~$l$ distinct primes. Alternatively,
set $V_0=\bigotimes_{i=1}^k V(\lambda_i)^{\tensor m_i}$.
Then
\begin{equation}\label{CH10.2}
\dim L^s(V)_{\nu+r\delta}=\frac1m\sum_{d|m}
\varphi_{r-s}(d) \dim (V_0)^{\tensor m/d}_{\nu/d}.
\end{equation}

\subl{CH20}
Retain the notations of~\secref{QL} and define, for any~$d$ dividing~$m$,
$$
\bpi^{1/d}=(\pi_1^{1/d},\dots,\pi_n^{1/d}),
$$
where~$\pi^{1/d}_i(u)=\prod_{j=0}^{m/d-1} \pi^0_i(\zeta^{jd}u)$.
In particular, if $d=1$ then~$\bpi^1=\bpi$.
It follows immediately from~\thmref{tenspr} and~\lemref{factor} that
$V(\bpi^{1/d})$ is isomorphic to~$V(\bpi^0)\tensor
V(\bpi^0_{\zeta^d})\tensor\cdots\tensor V(\bpi^0_{\zeta^{m-d}})$. Then
we conjecture the following quantum analogue of~\eqref{CH10.1}
and~\eqref{CH10.2}
\begin{conj}
Let~$\nu\in P$, $r\in\bz$ and~$s=0,\dots,m-1$. Then
\begin{align*}
\dim L^s(V(\bpi))_{\nu+r\delta}&=\frac1m\,\sum_{d|m} \varphi_{r-s}(d)
\dim V(\bpi^{1/d})_{\nu/d}\\
&=
\frac1m\,\sum_{d|m} \varphi_{r-s}(d) \dim (V(\bpi^0)^{\tensor m/d})_{\nu/d}.
\end{align*}
\end{conj}
In the rest of this section we prove this conjecture in
some special case for~$\lie g\cong
\lie{sl}_{n+1}$.

\subl{CH30}
Recall that~$V(\bpi)$ is a direct sum of eigenspaces of~$\eta_{\bpi^0}$
and that the eigenvalues of the latter are $m$th roots of unity.
Let~$V(\bpi)^{(k)}\subset V(\bpi)$ be the eigenspace of~$\eta_{\bpi^0}$
corresponding to the eigenvalue~$\zeta^k$.
Define, for all $v\in V(\bpi)$ and for all~$s\in\bz$
$$
\Pi_s(v)=\frac 1m\,\sum_{k=0}^{m-1} \zeta^{-ks}
\eta_{\bpi^0}^k(v)
$$
Evidently, $\Pi_s$ depends only on the residue class of~$s\pmod m$.
Observe that the~$\Pi_s$ are $\bu_q(\lie g)$-module endomorphisms of~$V(\bpi)$.
Furthermore, for all~$v\in V(\bpi)$ and for all~$r,s\in\bz$,
define
$$
\widehat \Pi_s(v\tensor t^r)=\Pi_{s-r}(v)\tensor t^r.
$$
\begin{lem}\label{proj}
\begin{enumerit}
\item The~$\Pi_s\,:\,s=0,\dots,m-1$ are orthogonal projectors onto
$V(\bpi)^{(s)}$.
\item The~$\widehat \Pi_s\,:\,s=0,\dots,m-1$ are $\bu_q^e$-module
endomorphisms of~$L(V(\bpi))$ and orthogonal projectors onto its
simple $\bu_q^e$ submodules~$L^s(V(\bpi))$.
\end{enumerit}
\end{lem}
\begin{proof}
For the first part, we have, for all~$v\in V(\bpi)$,
$$
\eta_{\bpi^0}(\Pi_s(v))=\frac 1m\,\sum_{k=0}^{m-1} \zeta^{-ks}
\eta_{\bpi^0}^{k+1}(v)=\frac 1m\,\sum_{k=1}^m \zeta^{-(k-1)s}
\eta_{\bpi^0}^k(v)=\zeta^s \Pi_s(v),
$$
whence~$\Pi_s(v)\in V(\bpi)^{(s)}$. Furthermore, write~$u=\Pi_s(v)$. Then
$\eta_{\bpi^0}^k(u)=\zeta^{ks} u$, and so
$$
\Pi_{s'} \Pi_s(v)=\frac 1m \sum_{k=0}^{m-1} \zeta^{k(s-s')} u.
$$
If~$s-s'=0\pmod m$ then the sum in the right-hand side equals~$m$,
whence~$\Pi_s^2=\Pi_s$. Otherwise, the sum
equals~$(\zeta^{m(s-s')}-1)/(\zeta^{(s-s')}-1)=0$, that
is~$\Pi_{s'}\circ \Pi_s=0$ if~$s\not=s'\pmod m$.

For the second part, take~$x\in(\bu_q)_k$, $v\in V(\bpi)$, $r\in\Z$. Then
\begin{align*}
\widehat \Pi_s (x(v\tensor t^r))&=
\widehat \Pi_s (x.v\tensor t^{r+k})=
\frac 1m\, \sum_{j=0}^{m-1}
\zeta^{j(r+k-s)} \eta_{\bpi^0}^j(x.v)\tensor t^{r+k}
\\
&=\frac 1m\,\sum_{j=0}^{ m-1}
\zeta^{j(r-s)} x.\eta_{\bpi^0}^j(v)\tensor t^{r+k}=
x.\widehat \Pi_{s}(v\tensor t^r).
\end{align*}
It follows that~$\widehat \Pi_s$ commutes with the action of~$\bu_q^e$.
Furthermore,
\begin{align*}
\eta(\widehat \Pi_s(v\tensor t^r))&
=\zeta^r \eta_{\bpi^0} \Pi_{s-r}(v)\tensor t^r
=\zeta^s \Pi_{s-r}(v)\tensor t^r
\\
&=\zeta^s\widehat \Pi_s(v\tensor t^r),
\end{align*}
whence the image of~$\widehat \Pi_s$ is contained in~$L^s(V(\bpi))$.
Finally, 
$$
\widehat \Pi_s' \widehat \Pi_s(v\tensor t^r)=\Pi_{s'-r}\Pi_{s-r}(v)
\tensor t^r.
$$
It follows immediately from the first part that~$\widehat
\Pi_s^2=\widehat \Pi_s$ whilst~$\widehat \Pi_{s'}\circ \widehat \Pi_s=0$ if~$s\not=s'
\pmod m$.
\end{proof}
Observe that~$\eta_{\bpi^0}$ preserves weight spaces of~$V(\bpi)$.
\begin{cor}
Given~$\nu\in P$, set~$V(\bpi)_\nu^{(k)}:=V(\bpi)_\nu\cap
V(\bpi)^{(k)}$. Then
$$
\dim L^s(V(\bpi))_{\nu+r\delta}=\dim V(\bpi)^{(k)}_\nu,
$$
where~$k=s-r\pmod m$.
\end{cor}
\begin{pf}
This follows immediately from~\lemref{proj}.
\end{pf}

\subl{CH40}
Throughout the rest of this section $\frak{g}$ is assumed to be
isomorphic to~$\frak{sl}_{n+1}$. It will be convenient to identify
the set~$\hat I$ with~$\bz/(n+1)\bz$ in the sense
that~$i+k\,:\, i\in\hat I$, $k\in\bz$ is understood as~$i+k\pmod{n+1}$.
Let~$V$ be the quantum analogue of the natural
representation of~$\frak g$. Explicitly, $V$ is an~$(n+1)$-dimensional
vector space over~$\bc(q)$ with a basis~$v_0,\dots,v_n$, the action
of the generators~$E_i, F_i, K_i\,:\,i\in I$ of~$\bu_q(\lie g)$
being given by
$$
E_i v_j=\delta_{i,j} v_{j-1},\quad F_i v_j=\delta_{i,j+1} v_{j+1},
\quad K_i v_j=q^{\delta_{j+1,i}-\delta_{j,i}} v_j,\qquad
i\in I,\quad j\in\hat I.
$$
Given~$a\in \bc(q)^\times$, we can endow~$V$ with a structure of
a~$\bu_q$-module which we denote by~$V(a)$ by setting
$$
E_0 v_j = \delta_{j,0} a v_n,\qquad F_0 v_j = \delta_{j,n} a^{-1}\,v_0,
\qquad K_0 v_j = q^{\delta_{j+1,0}-\delta_{j,0}} v_j,\qquad j\in\hat I.
$$
This module corresponds to the
$n$-tuple of polynomials~$\bpi=(1-au,1,\dots,1)$.

Let~$a,b\in\bc(q)^\times$ such that~$b/a\not=q^\bz$. Then, by
\cite{C}, $V(a)\tensor V(b)$ is irreducible as a~$\bu_q$-module
and~$V(a)\tensor V(b)\cong V(b)\tensor V(a)$.
A computation analogous to that of~\cite[5.4]{CPqa} and
based on a simple observation what~$V^{\tensor 2}\cong V(2\varpi_1)\oplus
V(\varpi_2)$ (where~$\varpi_2=0$ if~$n=1$) as a $\bu_q(\lie g)$-module,
allows one to obtain explicit formulae for the isomorphism~$I_z:
V(a)\tensor V(b)\to V(b)\tensor V(a)$ preserving the highest weight vector
\begin{align}
&I_z(v_i\tensor v_i)=v_i\tensor v_i,\qquad i\in\hat I\nonumber\\
&I_z(v_i\tensor v_j)=
\left(\frac{1-q^2}{z-q^2}\right) v_i\tensor v_j+
\left(\frac{q(z-1)}{z-q^2}\right) v_j\tensor v_i\label{tau}\\
&I_z(v_j\tensor v_i)=
\left(\frac{q(z-1)}{z-q^2}\right)v_i\tensor v_j+
\left(\frac{(1-q^2)z}{z-q^2}\right) v_j\tensor v_i,\qquad
0\le i<j\le n.\nonumber
\end{align}
where~$z=b/a$. In particular, we may always assume, without loss
of generality, that~$a=1$.
More generally, take~$a_1,\dots,a_m\in\bc(q)^\times$
such that~$a_i/a_j\not=q^\bz$. Then
$I_{a_{i+1}/a_i,i}:=\id^{\tensor i-1}\tensor I_{a_{i+1}/a_i}
\tensor \id^{\tensor m-i-1}$
gives an isomorphism of~$\bu_q$-modules with permutes
the~$i$th and the $(i+1)$th factors in the tensor product
$V(a_1)\tensor\cdots\tensor V(a_m)$. The latter is isomorphic
to~$V(\bpi)$ where~$\bpi=( \prod_{i=1}^s (1-a_i u),1,\dots,1)$.

Now, suppose that~$a_i=\zeta^{i-1}$, where~$\zeta$ is an~$m$th
primitive root of unity. Then we are in the situation of~\lemref{factor}
and~$\bpi^0=(1-u,1,\dots,1)$. It follows that the isomorphism
$\tau_{\bpi^0}$ can be constructed explicitly
as~$\mathbf I_m=
I_{\zeta^{m-1},m-1}\circ I_{\zeta^{m-2},m-2}\circ\cdots
\circ I_{\zeta,1}$. Furthermore, since~$V(\zeta^i)$ is simple
as a~$\bu_q(\lie g)$-module,
$\phi_{\bpi^0_{\zeta^i},\zeta^{-1}}=\id$ as a map of vector spaces.
Therefore, $\mathbf I_m$ can be identified with the map~$\eta_{\bpi^0}$
defined in \ref{QL50}.

\subl{CH50}
Recall that~$\bc(q)$ is the quotient field of~$\bc[q]$ and define
a subring~$A$ of~$\bc(q)$ by~$A:=\{f/g\,:\,
f,g\in\bc[q],\, g(0)\not=0\}$. Then~$A$
is a local ring, the unique maximal ideal being~$qA$. Evidently,
if~$f/g\in A$, then~$f/g=f(0)/g(0)\pmod{qA}$. In particular, $A/qA\cong\bc$.
Observe also that, given an~$A$-module~$M$, we get~$\overline{M}:=
M/qM\cong M\tensor_A
A/qA$. Given $m\in M$, we denote its canonical
image in~$\overline{M}$ by~$\bar m$.

Let~$\bb(\bpi)$ be the set of all tensor products
of the form~$v_{i_m}\tensor\cdots\tensor v_{i_1}\,:\,i_r\in\hat I$
and set~$\cal L(\bpi)=\bigoplus_{\bov\in\bb(\bpi)} A\bov$.
It follows from~\eqref{tau} that
$\eta_{\bpi^0}$ maps~$\cal L(\bpi)$ into itself and in particular
preserves~$q\cal L(\bpi)$. Furthermore, if~$\nu\in P$ is a weight of~$V(\bpi)$,
set~$\cal L(\bpi)_\nu:=\cal L(\bpi)\cap V(\bpi)_\nu$. Since
each~$\bov\in\bb(\bpi)$
is a weight vector, it follows that $\eta_{\bpi^0}$
preserves~$\cal L(\bpi)_\nu$. Set~$\overline{\cal L(\bpi)}_\nu:=
\overline{\cal L(\bpi)_\nu}$. Then~$\dim_\bc \overline{\cal L(\bpi)}_\nu=
\operatorname{rank}_A \cal L(\bpi)_\nu=\dim_{\bc(q)}V(\bpi)_\nu$.

\begin{lem}\label{maj}
Given~$\bov=v_{i_m}\tensor\cdots\tensor v_{i_1}\in \bb(\bpi)$, define
$$
\desc(\bov)=\{r\,:\, i_{r+1}<i_r,\,
1\le r< m\}
$$
and let~$\Maj(\bov)$ be the sum of elements
of~$\desc(\bov)$ if the latter is non-empty and
zero otherwise. Then
$$
\eta_{\bpi^0}(\bov)=\zeta^{\Maj(\bov)} \bov\pmod{q\cal L(\bpi)}.
$$
\end{lem}
\begin{pf}
Observe that~\eqref{tau} yields, for all~$r=1,\dots,m-1$,
$$
I_{\zeta^{m-r},m-r}(\bov)=
\begin{cases}
\bov,\qquad& i_{r+1}\ge i_r\\
\zeta^{r-m} \bov,&
i_{r+1} < i_r
\end{cases}
\pmod{q\cal L(\bpi)},
$$
It follows immediately from the definition of~$\desc(\bov)$ that
$I_{\zeta^{m-r},m-r}(\bov)=\zeta^r \bov\pmod{q\cal L(\bpi)}$
provided that~$r\in\desc(\bov)$
and~$I_{\zeta^{m-r},m-r}(\bov)=\bov
\pmod{q\cal L(\bpi)}$ otherwise. The result is now immediate
since~$\eta_{\bpi^0}(q\cal L(\bpi))\subset q\cal L(\bpi)$.
\end{pf}

\subl{CH60}
Denote by~$\bar\eta_{\bpi^0}$ the map of $\bc$-vector space
$\overline{\cal L(\bpi)}$ into itself which obtains canonically
from~$\eta_{\bpi^0}$. By the above Lemma, the
eigenvalues of~$\bar\eta_{\bpi^0}$ are $m$th roots of unity.
\begin{lem}\label{dim}
Let~$\overline{\cal L(\bpi)}\,{}^{(k)}$
be the eigenspace of~$\bar\eta_{\bpi^0}$
corresponding to the eigenvalue~$\zeta^k$. Then, for all~$k=0,\dots,m-1$,
$$
\dim_{\bc(q)} V(\bpi)^{(k)}=\dim_\bc \overline{\cal L(\bpi)}\,{}^{(k)}.
$$
Moreover, if we set~$\overline{\cal L(\bpi)}\,{}^{(k)}_\nu:=
\overline{\cal L(\bpi)}\,{}^{(k)}\cap \overline{\cal L(\bpi)}_\nu$,
then
$$
\dim_{\bc(q)} V(\bpi)^{(k)}_\nu=\dim_\bc \overline{\cal L(\bpi)}\,{}^{(k)}_\nu.
$$
\end{lem}
\begin{pf}
Given~$\bov\in\bb(\bpi)$, set
$\langle \bov\rangle:=\Pi_k(\bov)$ where~$k=\Maj(\bov)$.
Evidently, $\langle\bov\rangle\in \cal L(\bpi)$ and, moreover,
$\langle
\bov\rangle=\bov\pmod{q\cal L(\bpi)}$ by~\lemref{maj}.
Since the canonical
images of~$\bov\in\bb(\bpi)$ in~$\overline{\cal L(\bpi)}$ form a basis
of $\overline{\cal L(\bpi)}$ over~$\bc$,
the $\langle\bov\rangle\,:\, \mathbf
v\in\bb(\bpi)$ generate~$\cal L(\bpi)$ as an~$A$-module and
are linearly independent over~$A$ by Nakayama's Lemma.
Therefore, the $\langle\bov\rangle\,:\, \mathbf
v\in\bb(\bpi)$ form a basis of~$V(\bpi)$ over~$\bc(q)$.
Furthermore, by~\lemref{proj}(i), $\langle\bov\rangle\in V(\bpi)^{(k)}$.
We conclude that
$V(\bpi)^{(k)}$ contains a linearly independent subset
$\{\langle \bov\rangle\,:\, \bov\in\bb(\bpi),\,
\Maj(\bov)=k\pmod m\}$, whose
cardinality equals~$\dim_\bc\overline{\cal L(\bpi)}\,{}^{(k)}$ by~\lemref{maj}.
Thus, $\dim_{\bc(q)} V(\bpi)^{(k)}\ge \dim_\bc
\overline{\cal L(\bpi)}\,{}^{(k)}$. On the other hand, since $V(\bpi)$
(respectively, $\overline{\mathcal L(\bpi)}$) is a direct sum of eigenspaces
of $\eta_{\bpi^0}$ (respectively, $\bar\eta_{\bpi^0}$), it follows that
$\dim_{\bc(q)}V(\bpi)=\sum_{k=0}^{m-1} \dim_{\bc(q)} V(\bpi)^{(k)}
\ge \sum_{k=0}^{m-1} \dim_\bc\overline{\cal L(\bpi)}\,{}^{(k)}=
\dim_{\bc}\overline{\cal L(\bpi)}=\dim_{\bc(q)}V(\bpi)$, whence
the desired equality. The second assertion is immediate since
all the~$\bov\in\bb(\bpi)$ are weight vectors.
\end{pf}

\subl{CH65}
It follows immediately from~\lemref{maj} that
$$
\dim_\bc \overline{\cal L(\bpi)}\,{}^{(k)}_\nu=\#\{\bov\,:\, \bov\in
V(\bpi)_\nu,\,\Maj(\bov)=k\pmod m\}.
$$
The cardinality of the set which appears in the right-hand side
was computed in~\cite{G}. Namely,
there is a bijection
between the set of weights of~$V(\bpi)$ (or~$\cal L(\bpi)$) and
the set
$$
\{(k_0,\dots,k_n)\in\bn^{n+1}\,:\,k_0+\cdots+k_n=m\}.
$$
Indeed, take $\bov=v_{i_m}\tensor\cdots\tensor v_{i_1}\in\bb(\bpi)$
and set $k_i(\bov)\,:\,i=0,\dots,n$
where~$k_i(\bov)=\#\{ r\,:\,
i_r=i\}$. Then~$\bov$ is of weight~$\nu=\sum_{i=0}^n
k_i(\varpi_{i+1}-\varpi_i)=
\sum_{i=1}^n (k_{i-1}-k_i)\varpi_i$, where we set~$\varpi_0=\varpi_{n+1}=0$.
It follows from the definition that~$k_0(\bov)+\cdots+k_i(\bov)=m$.
One can easily check that~$\bov,\bov'\in\bb(\bpi)$ are of the same
weight if and only if~$k_i(\bov)=k_i(\bov')$ for all~$i=0,\dots,n$.

Suppose that~$\nu$ corresponds to~$(k_0,\dots,k_n)$.
Then
\begin{align*}
&\dim_\bc \overline{\cal L(\bpi)}_\nu = \dim_{\bc(q)} V(\bpi)_\nu=
\binom{m}{k_0,\dots,k_n}
\\
\intertext{and by~\cite[A.1-A.6]{G}}
&\dim_\bc \overline{\cal L(\bpi)}\,{}^{(k)}_\nu=\frac1m \sum_{d|m}
\varphi_k(d)\binom{\frac md}{\frac{k_0}{d},\dots,\frac{k_n}{d}}=
\frac1m \sum_{d|m} \varphi_k(d) \dim_{\bc(q)} V(\bpi^{1/d})_{\nu/d}.
\end{align*}
Applying~\ref{CH30} and~\lemref{dim}, we obtain
the following
\begin{prop}
Assume that~$\lie g$ is isomorphic to~$\lie{sl}_{n+1}$.
Let~$\bpi^0=(1-u,1,\dots,1)$ and
$\bpi=(\prod_{j=0}^{m-1} (1-\zeta^{j-1}u),1,\dots,1)$, where~$\zeta$
is an~$m$th primitive root of unity. Then
\begin{align*}
\dim_{\bc(q)} L^s (V(\bpi))_{\nu+r\delta}&=
\frac1m \sum_{d|m} \varphi_{r-s} (d) \dim_{\bc(q)} V(\bpi^{1/d})_{\nu/d}\\
&=\frac1m\sum_{d|m} \varphi_{r-s}(d) \dim_{\bc(q)}
(V(\bpi^0)^{\tensor m/d})_{\nu/d},
\end{align*}
for all~$\nu\in P$, $r\in\bz$ and~$s=0,\dots,m-1$.
\end{prop}

\subl{CH70}
In the remainder of this section we will discuss the crystal basis theory
of our modules~$L(V(\bpi))$. One can easily see that~$V(\zeta^i)$ does not a
admit a crystal basis in the sense of Kashiwara. However, one
can slightly modify the definition of a crystal basis so that
it makes sense in our particular case.

Let~$M$ be an integrable~$\wh\bu_q$-module or
a finite dimensional~$\bu_q$-module. Then, for~$i\in\hat I$ fixed,
the Kashiwara operators~$\tilde E_i$, $\tilde F_i$ are defined
in the following way (cf.~\cite[2.2]{Ka1} or~\cite[16.1]{L}). 
Any~$m\in M_\lambda$ can be written uniquely as~$m
=\sum_{s\ge0, s+t\ge0} F_i^{(s)}u_s$, where~$F_i^{(s)}:=F_i^s/[s]!$,
$t=(\lambda\,|\,\alpha_i)$, $u_s\in \ker E_i\cap M_{\lambda+s\alpha_i}$
and equals zero for~$s$ sufficiently large.
Then
$$
\tilde E_i m = \sum_{s>0,\,s+t\ge0} F_i^{(s-1)}u_s,\qquad
\tilde F_i m = \sum_{s\ge0,\,s+t\ge0} F_i^{(s+1)}u_s.
$$
\begin{defn}
Let~$\mathcal L$ be a free $A$-submodule of~$M$ satisfying~$M=\cal L\tensor_A
\bc(q)$ and~$\cal B$ be a basis
of~$\overline{\cal L}$. Fix~$\zeta\in\bc^\times$.
We say that~$(\cal L,\cal B)$ is a $\zeta$-crystal
basis of~$M$ if
\begin{enumerit}
\item $\cal L=\bigoplus_\nu \cal L_\nu$, $\cal B=
\coprod_\nu \cal B_\nu$, where~$\cal L_\nu=\cal L\cap M_\nu$ 
and~$\cal B_\nu=\cal B\cap\overline{\cal L}_\nu$.
\item $\cal L$ is stable by the~$\tilde E_i$, $\tilde F_i$
for all~$i\in\hat I$, hence the~$\tilde E_i$, $\tilde F_i$ act on~$\overline{
\cal L}$.
\item $\tilde E_i\cal B,\tilde F_i\cal B
\subset \zeta^{\bz \delta_{i,0}}\cal B\cup\{0\}$, 
for all~$i\in\hat I$.
\item For $\bov,\bov'\in\cal B$ one has
$\bov'=\zeta^k \tilde F_i\bov$ if and only if~$\bov=\zeta^{-k}\tilde E_i
\bov'$.
\end{enumerit}
\end{defn}
If~$\zeta=1$, the above definition reduces to the Kashiwara's definition of
crystal bases (cf.~\cite[2.3]{Ka1}).

The following Lemma is an adaptation of~\cite[Lemma~20.1.2]{L}
for~$\zeta$-crystal bases.
\begin{lem}
Assume that~$(\cal L,\cal B)$ is a $\zeta$-crystal
basis of~$M$.
Let~$m\in\cal L_\lambda$ and fix~$i\in\hat I$. Write~$m=\sum_{s\ge0\,
s+t\ge0} F_i^{(s)} u_s$, where~$u_s\in\ker E_i\cap M_{\lambda+s\alpha_i}$,
and~$u_s=0$ if~$s\gg0$.
Then
\begin{enumerit}
\item For all~$s\ge0$ and~$r\ge0$, $F_i^{(r)}x_s\in\cal L$.
\item If~$\bar m\in\cal B$, then there exists~$s_0$ such that~$u_s\in
q\cal L$ if~$s_0\not=s$, $\bar u_{s_0}\in\zeta^{\bz\delta_{i,0}}
\cal B$ and~$m=F_i^{(s_0)} u_{s_0}\pmod{q\cal L}$.
\end{enumerit}
\end{lem}
\begin{pf}
The argument is
an obvious modification of that of~\cite[Lemma~20.1.2]{L}.
\end{pf}

\subl{CH80}
Retain the notations of~\ref{CH40}.
\begin{prop}
Fix~$\zeta\in\bc^\times$ and let~$\bpi=(\prod_{j=0}^{N-1}(1-\zeta^j u),1,
\dots,1)$. Let~$V(\bpi)$ be the corresponding simple $\bu_q$-module
which is isomorphic to~$V(1)\tensor\cdots\tensor V(\zeta^{N-1})$.
Set~$\bb(\bpi)=\{ v_{i_1}\tensor\cdots\tensor v_{i_N}\,:\, i_r\in\hat I\}$ 
and~$\cal L(\bpi)=
\bigoplus_{\bov\in \bb(\bpi)} A \bov$. Then
the pair~$(\mathcal L(\bpi), \bb(\bpi))$ 
forms a $\zeta$-crystal basis of~$V(\bpi)$. Moreover, the
action of~$\tilde E_i$, $\tilde F_i$ on~$\bov=v_{i_1}\tensor
\cdots\tensor v_{i_N}$ is given by
\begin{align*}
\tilde E_i\bov=\zeta^{(r-1)\delta_{i,0}} &v_{i_1}
\tensor\cdots\tensor \tilde E_i v_{i_r}\tensor\cdots\tensor v_{i_N}
\pmod{q\cal L(\bpi)},
\\
\tilde F_i\bov=\zeta^{-(s-1)\delta_{i,0}} &v_{i_1}
\tensor\cdots\tensor \tilde F_i v_{i_s}\,\tensor\cdots\tensor v_{i_N}
\pmod{q\cal L(\bpi)},
\end{align*}
where~$r$ and~$s$ are determined by the standard Kashiwara rules for
the tensor product of crystals {\em({\em cf.~\cite[Theorem~1]{Ka1}}
and~\cite[1.3]{Ka1a})}.
\end{prop}
\begin{pf}
The property~{(i)} of~\defref{CH70} is obvious. The other three
for~$i\in I$ follow from the standard results of Kashiwara on crystal bases
(cf. for example~\cite[20.2]{L}). 
So, it only remains to prove~{(ii)-(iv)} for~$i=0$. 

The proof is basically an adaptation of the standard argument. Let~$\bu_0$
be the subalgebra of~$\wh\bu_q$ generated by~$E_0$, $F_0$ and~$K_0^{\pm1}$,
which is isomorphic to~$\bu_q(\lie{sl}_2)$. Throughout the
rest of the proof we shall omit indices of the operators~$E_0$ and~$F_0$.

We use the following inductive argument. Set~$M_k:=V(1)\tensor
V(\zeta)\tensor\cdots\tensor V(\zeta^{k-1})$. 
It is clear that~\hbox{(ii)--(iv)} of~\defref{CH70}
hold for~$M_1$ and that~{(iv)} follows from~{(ii)--(iii)}, the second
assertion of our Proposition and Kashiwara's tensor product rules.
Now, suppose that they hold for~$M_k$
and the 
pair~$(\mathcal L_k,\bb_k)$,
where~$\bb_k=\{v_{i_1}\tensor \cdots\tensor v_{i_k}\,:\, i_r\in\hat I\}$ and
$\mathcal L_k=\bigoplus_{\bob\in\bb_k} A\bob$. 
Given~$\bob\in\bb_k$ of weight~$\nu$,
set~$t=(\nu\,|\,\alpha_0)$ and write~$\bob=\sum_{s\ge0,s+t\ge0} 
F^{(s)} u_s$ as in~\ref{CH70}.

Notice that~$M_{k+1}\cong V(1)\tensor \phi_{\bpi_k,\zeta}M_k$,
where~$\bpi_k=(\prod_{j=0}^{k-1} (1-\zeta^j u),1,\dots,1)$. 
For the inductive step we should prove first that, for all~$i\in\hat I$ and
for all~$\bob\in\bb_k$,
$\tilde E (v_i\tensor\bob),\tilde F(v_i\tensor\bob)\in\cal L_{k+1}$.
Since the~$v_i\,:\,i\not=0,n$ span trivial~$\bu_0$-modules, we
can write
$$
v_i\tensor\bob=\sum_{s\ge0,s+t\ge0} v_i\tensor F^{(s)}u_s
=\sum_{s\ge0,s+t\ge0} F^{(s)}(v_i\tensor \zeta^s u_s).
$$
Since~$E(v_i\tensor u_s)=0$ if~$i\not=0$,
the claim follows by induction hypothesis. When it is easy to check
that~{(iii)} holds for~$\bov\in\bb_{k+1}$ of the form~$v_i\tensor
\bob$ where~$\bob\in\bb_k$ and~$i\not=0,n$. Now, by~\lemref{CH70},
$\bob=F^{(s_0)}u_{s_0}\pmod{q\cal L_k}$ for some~$s_0\ge0$. It follows
that~$v_i\tensor\bob=\zeta^{s_0} F^{(s_0)}(v_i\tensor u_{s_0})$. If~$s_0$=0,
then obviously~$\tilde E(v_i\tensor\bob)=0$. Otherwise,
$\tilde E(v_i\tensor\bob)=\zeta^{s_0} F^{(s_0-1)}(v_i\tensor u_{s_0})
\pmod{q\cal L_{k+1}}=\zeta v_i\tensor \tilde E\bob\pmod{q\cal L_{k+1}}$.
Similarly, $\tilde F(v_i\tensor\bob)=\zeta^{s_0} F^{(s_0+1)}(v_i\tensor
\bob)\pmod{q\cal L_{k+1}}=\zeta^{-1} 
v_i\tensor\tilde F\bob\pmod{q\cal L_{k+1}}$, which proves 
the second assertion and~{(iv)}. 

It remains to consider the case when~$i=0,n$. 
One can easily check that the weight vectors
$$
X_s:=v_n\tensor u_s,\qquad
Y_s:=v_n\tensor Fu_s-\zeta q^{t+2s}[t+2s] v_0\tensor u_s
$$
are annihilated by~$E$. 
Furthermore, for all~$r>0$,
\begin{align*}
F^{(r)} X_s &= \zeta^{-r} q^r v_n\tensor F^{(r)}u_s+
\zeta^{-r+1} v_0\tensor F^{(r-1)}u_s,
\\
F^{(r)} Y_s &= \zeta^{-r} q^r [r+1] v_n \tensor F^{(r+1)}u_s+
\zeta^{-r+1}([r]-q^{t+2s-r}[t+2s]) v_0 \tensor F^{(r)}u_s
\\
&=\zeta^{-r}\, \frac{q^{2r}-1}{q^2-1}\,v_n\tensor F^{(r+1)}u_s+
\zeta^{-r+1} q^r\,\frac{q^{2(t+2s-r)}-1}{q^2-1}\, v_0\tensor F^{(r)}u_s.
\\
\intertext{In particular, we have}
F^{(s)} X_s &= \zeta^{-s} q^s v_n\tensor F^{(s)}u_s+
\zeta^{-s+1} v_0\tensor F^{(s-1)}u_s,
\\
F^{(s-1)} Y_s 
&=\zeta^{-s+1}\,\frac{q^{2(s-1)}-1}{q^2-1}\, v_n\tensor F^{(s)} u_s+
\zeta^{-s+2} q^s \,\frac{q^{2(t+s+1)}-1}{q^2-1}\, v_0 \tensor F^{(s-1)}u_s.
\end{align*}
Fix~$s>0$. Evidently, $v_n\tensor F^{(s)} u_s$ and~$v_0 \tensor F^{(s-1)}u_s$
are non-zero and linearly independent. Furthermore, the determinant of the
matrix of coefficients of~$F^{(s)} X_s$ and~$F^{(s-1)} Y_s$ in
the basis consisting of~$v_n\tensor F^{(s)} u_s$ and~$v_0 
\tensor F^{(s-1)}u_s$ equals~$-\zeta^{-2(s-1)}\pmod{qA}$, 
whence is a unit in~$A$.
Therefore, there exist~$a_s, b_s \in A$ such that~$v_n\tensor F^{(s)} u_s=
a_s F^{(s)} X_s+b_s F^{(s-1)} Y_s$, $s>0$. Thus, we can write
\begin{align*}
v_n\tensor\bob&=\sum_{s\ge0,\,s+t\ge0} v_n\tensor F^{(s)}u_s=
\sum_{s\ge0,\,s+t\ge0} a_s F^{(s)} X_s+\sum_{s>0,\,
s+t\ge0} b_s F^{(s-1)} Y_s\\&=
\sum_{s\ge0,\,s+t+1\ge0} F^{(s)} (a_s X_s+b_{s+1} Y_{s+1})=
\sum_{s\ge0,\,s+t+1\ge0} F^{(s)} U_s,
\end{align*}
which is the decomposition of~$v_n\tensor\bob$ of~\ref{CH70}.
It follows from the definitions of~$X_s$, $Y_s$,
\lemref{CH70} and the induction hypothesis that~$\tilde E(v_n\tensor\bob)=
\sum_{s>0,\,s+t+1\ge0} F^{(s-1)}U_s$ and~$\tilde F(v_n\tensor\bob)
=\sum_{s\ge0,s+t+1\ge0} F^{(s+1)}U_s$ lie in~$\mathcal L_{k+1}$.
A similar argument shows that~$\tilde E(v_0\tensor\bob),
\tilde F(v_0\tensor\bob)\in\mathcal L_{k+1}$.

Now observe that~$F^{(s)}X_s=\zeta^{-s+1}v_0\tensor F^{(s-1)}u_s
\pmod{q\cal L_{k+1}}$ whilst~$F^{(s-1)}Y_s=
\zeta^{-s+1} v_n\tensor F^{(s)} u_s\pmod{q\cal L_{k+1}}$, $s>0$.
Assume first that~$\tilde E\bob=0\pmod{q\cal L_k}$. Then~$\bob=u_0
\pmod{q\cal L_k}$ and so~$v_n\tensor\bob=X_0\pmod{q\cal L_{k+1}}$.
It follows that~$\tilde E(v_n\tensor\bob)=0=\tilde Ev_n\tensor \bob
\pmod{q\cal L_{k+1}}$. On the other hand, $\tilde F(v_n\tensor\bob)=
v_0\tensor u_s\pmod{q\cal L_{k+1}}=\tilde Fv_n\tensor\bob$. Suppose
further that~$\tilde E\bob\notin q\cal L_k$. Then~$\bob=F^{(s)}u_{s}
\pmod{q\cal L_k}$ for some~$s>0$. It follows that~$v_n\tensor
\bob=\zeta^{s-1} F^{(s-1)}Y_s\pmod{q\cal L_{k+1}}$. If~$s=1$, then
$\tilde E(v_n\tensor\bob)=0=\tilde Ev_n\tensor\bob\pmod{q\cal L_{k+1}}$.
Otherwise, 
$$
\tilde E(v_n\tensor\bob)=\zeta^{s-1} F^{(s-2)}Y_s
=\zeta v_n\tensor F^{(s-1)}u_s=
\zeta v_n\tensor\tilde E\bob\pmod{q\cal L_{k+1}}.
$$
Similarly, $\tilde F(v_n\tensor\bob)=\zeta^{s-1} F^{(s)}Y_s=
\zeta^{-1} v_n\tensor\tilde F\bob\pmod{q\cal L_{k+1}}$.
We omit an analogous computation for~$v_0\tensor\bob$.
\end{pf}

\subl{CH90}
Suppose now that~$\zeta$ is an $m$th primitive root of unity 
and retain the notations of~\ref{CH50}--\ref{CH60}.
\begin{prop}
Set~$\bb(\bpi)^{(k)}=\{\bov\in\bb(\bpi)\,:\, \Maj(\bov)=k\pmod m\}$
and define
$$
\wh{\cal L^s(\bpi)}=\bigoplus_{r\in\bz}\mskip20mu\bigoplus_{\bov\in
\bb(\bpi)^{(s-r)}} A\langle\bov\rangle\tensor t^r,
\qquad \wh{\bb^s(\bpi)}=\coprod_{r\in\Z} \bb(\bpi)^{(s-r)}\tensor t^n.
$$
Then~$(\wh{\cal L^s(\bpi)}, \wh{\bb^s(\bpi)})$ is a $\zeta$-crystal basis
of~$L^s(V(\bpi))$. 
\end{prop}
\begin{pf}
Set $\wh{\mathcal L}=\mathcal L(\bpi)\tensor_A A[t,t^{-1}]$,
$\wh{\mathcal B}=\coprod_{r\in\Z} \bb(\bpi)\tensor t^r$. One
can easily check that
the above proposition implies that 
$(\wh{\mathcal L},\wh{\mathcal B})$ 
is a $\zeta$-crystal basis of~$L(V(\bpi))$. 
Set~$\wh{\cal L^s}=\wh{\cal L}\cap L^s(V(\bpi))$. We claim that~$\wh{\cal L^s}=
\wh{\cal L^s(\bpi)}$ defined above. Indeed, since~$\langle\bov\rangle$
lies in~$\cal L(\bpi)$ for all~$\bov\in\bb(\bpi)$, 
$\wh{\cal L^s(\bpi)}\subset\wh{\cal L}$.
Suppose further that~$\bov\in
\bb(\bpi)^{(s-r)}$. Then~$\eta(\langle \bov\rangle\tensor t^r)=
\zeta^r \eta_{\bpi^0}(\langle\bov\rangle)\tensor t^r=
\zeta^s\langle \bov\rangle\tensor t^r$ (cf. the proof of~\lemref{CH60}), 
whence~$\wh{\cal L^s(\bpi)}
\subset L^s(V(\bpi))\cap \wh{\cal L}=\wh{\cal L^s}$ as an~$A$ submodule. 
Furthermore, both~$\wh{\cal L^s(\bpi)}$ and~$\wh{\cal L^s}$ are
direct sums of free weight $A$-submodules
and~$\wh{\cal L^s(\bpi)}_{\nu+r\delta}\subset \wh{\cal L^s}_{\nu+r\delta}$
for all~$\nu\in P$, $r\in\Z$.
Observe that~$\wh{\cal L^s(\bpi)}_{\nu+r\delta}$ is generated 
as an~$A$-module by~$\langle
\bov\rangle\tensor t^r$ where~$\bov\in\bb(\bpi)_\nu^{(s-r)}$. Since
the images of these elements in~$\wh{\cal L^s}_{\nu+r\delta}/q\wh{\cal L^s}
_{\nu+r\delta}$ form a basis
of that vector space by~\ref{CH60}, it follows by Nakayama's lemma
that they generate~$\wh{\cal L^s}$ as an~$A$-module. Therefore,
$\wh{\cal L^s(\bpi)}=\wh{\cal L^s}$. 
The result now follows immediately from~\propref{CH80}
and~\cite[Proposition~3.6 and Lemma~A.1]{G}.
\end{pf}

\newpage
\bibliographystyle{amsplain}

\end{document}